\documentclass[final,3p,times]{elsarticle}

\usepackage{algorithm}
\usepackage{algorithmicx}
\usepackage[noend]{algpseudocode}

\usepackage{amsmath}
\usepackage{amsfonts,amssymb}
\usepackage{calc}

\usepackage{subfig}
\usepackage{textcomp}
\usepackage{xspace}

\usepackage{graphicx}
\usepackage[pdftex]{hyperref}

\newcommand{\K}{\mathbb{K}}
\newcommand{\cpu}{\textsc{cpu}}
\newcommand{\gpu}{\textsc{gpu}}

\newcommand{\bem}{\textsc{bem}\xspace}
\newcommand{\fmm}{\textsc{fmm}\xspace}

\newcommand{\cpp}{C$^{++}$}

\newcommand{\bigO}{\mathcal{O}}
\renewcommand{\O}[1]{\mathcal{O}(#1)}

\newcommand{\ptom}{\textsc{p}\texttwooldstyle\textsc{m}\xspace} 
\newcommand{\ltop}{\textsc{l}\texttwooldstyle\textsc{p}\xspace} 
\newcommand{\mtop}{\textsc{m}\texttwooldstyle\textsc{p}\xspace} 
\newcommand{\mtom}{\textsc{m}\texttwooldstyle\textsc{m}\xspace} 
\newcommand{\mtol}{\textsc{m}\texttwooldstyle\textsc{l}\xspace} 
\newcommand{\ltol}{\textsc{l}\texttwooldstyle\textsc{l}\xspace}  
\newcommand{\ptop}{\textsc{p}\texttwooldstyle\textsc{p}\xspace} 

\newcommand{\ncrit}{N_{\text{CRIT}}}
\newcommand{\pmin}{p_{\text{min}}}
\newcommand{\tsolve}{t_{\text{solve}}}
\newcommand{\mac}{\textsc{mac}}

\newcommand{\gmres}{\textsc{gmres}\xspace}

\newcommand{\di}[1]{\text{d}#1}
\newcommand{\partiald}[2]{\frac{\partial #1}{\partial #2}}
\newcommand{\partialdi}[2]{\partial #1 / \partial #2}
\newcommand{\nhat}{\hat{n}}
\newcommand{\vect}[1]{\mathbf{#1}}

\graphicspath{{figs/}} 

\begin{document}

\begin{frontmatter}



\title{Inexact Krylov iterations and relaxation strategies with fast-multipole boundary element method\tnoteref{t1}}
\tnotetext[t1]{This work was partially supported by the National Science Foundation under award ACI-1149784.}

\author[gwu]{Tingyu Wang}
\ead{twang66@gwu.edu}

\author[bu,nvidia]{Simon K. Layton\fnref{fn1}}
\ead{slayton@nvidia.com}

\author[gwu]{Lorena A. Barba}
\ead{labarba@gwu.edu}

\fntext[fn1]{Currently at Nvidia, Corp.}

\address[gwu]{Department of Mechanical and Aerospace Engineering, The George Washington University, Washington DC, 20052}
\address[bu]{Department of Mechanical Engineering, Boston University, Boston, MA, 02215}
\address[nvidia]{Nvidia, Corp., Santa Clara, CA}

\begin{abstract}
Boundary element methods produce dense linear systems that can be accelerated via multipole expansions. Solved with Krylov methods, this implies computing the matrix-vector products within each iteration with some error, at an accuracy controlled by the order of the expansion, $p$. We take advantage of a unique property of Krylov iterations that allow lower accuracy of the matrix-vector products as convergence proceeds, and propose a relaxation strategy based on progressively decreasing $p$. In extensive numerical tests of the relaxed Krylov iterations, we obtained speed-ups of between $2.1\times$ and $3.3\times$ for Laplace problems and between $1.7\times$ and $4.0\times$ for Stokes problems. We include an application to Stokes flow around red blood cells, computing with up to 64 cells and problem size up to 131k boundary elements and nearly 400k unknowns. The study was done with an in-house multi-threaded C++ code, on a hexa-core CPU. The code is available on its version-control repository, \href{https://github.com/barbagroup/fmm-bem-relaxed}{https://github.com/barbagroup/fmm-bem-relaxed}.
\end{abstract}


\end{frontmatter}



\section{Introduction}

The boundary element method (\bem) is popular in applied mechanics for solving problems governed by the equations of Laplace, Stokes, Helmholtz and Lam{\'e}. Applications include electrostatics, low-Reynolds-number flow, acoustics and linear elasticity and span all scales from biomolecules, micro-electromechanical and biomedical devices to wind turbines, submarines and  problems in aerospace and geodynamics. The key feature of \bem\ is formulating the governing partial differential equations in equivalent boundary integral equations, and discretizing over the boundaries. This process reduces the dimensionality of the problem (three-dimensional problems are solved on two-dimensional surfaces), but generates dense systems of algebraic equations. 
The computational complexity of dense solvers, scaling as $\O{N^3}$ for direct methods and $\O{N^2}$ for iterative methods, frustrated large-scale applications of \bem\ until the late 1990s, when \bem\ researchers began working out how to incorporate fast algorithms 
\cite{Liu2006}. Treecodes and fast multipole methods (\fmm) reduce the complexity of \bem\ solutions to $\O{N \log N}$ and $\O{N}$, respectively, although often with serious programming effort. Nevertheless, large-scale \bem\ simulations are now possible, especially given the parallel scalability of the \fmm \cite{YokotaETal2011a,YokotaBarba2011a}.

Accelerating \bem\ solutions with \fmm\ hinges on seeing the dense matrix-vector products done within each iteration of the linear solver, as $N$-body problems. Gauss quadrature points on source panels and collocation points on target panels interact via the Green's function of the governing equation, similar to how charges, particles or masses interact under electrostatic or gravitational potentials. In the same way, far-away sources can be clustered together and represented by series expansions to calculate their contribution on a target point with controllable accuracy. As a result, the matrix-vector products are computed with some error. To ensure convergence of the iterative solver, or achieve a desired accuracy in the solution, we might require the order of the series expansions, $p$, to be sufficiently large. But Krylov methods have the surprising property of only requiring high accuracy on the first iteration, while accuracy requirements can be relaxed on later iterations. This property offers the opportunity to use different values of $p$ in the \fmm\ as the iterations progress to convergence, reducing the total time to solution.

This paper presents a relaxation strategy for fast-multipole boundary element methods consisting of decreasing the orders of expansion $p$ as Krylov iterations progress in the solver. We tested extensively using the Laplace and Stokes equations, and also with an application of Stokes flow around red blood cells. The aim of the study was to show that a \bem\ solution with inexact Krylov iterations converges, despite low-accuracy matrix-vector multiplications at later iterations, and to find out the speed-ups that can be obtained. We wrote an in-house  multi-threaded code in \cpp\ that enable experimenting in a variety of scenarios.\footnote{The code is available for replication of our results at \href{https://github.com/barbagroup/fmm-bem-relaxed}{https://github.com/barbagroup/fmm-bem-relaxed}.}

\section{Methods for the integral solution of elliptic equations using inexact {\small GMRES}}

\subsection{Boundary-integral solution of the Laplace equation}

To write the Laplace equation, $\nabla^{2}\phi(\vect{x}) = 0$,  in its integral formulation, we use the classical procedure of multiplying by the Green's function and integrating, applying the divergence theorem of Gauss and the chain rule, then dealing with singularities by a limiting process. This results in
\begin{equation}\label{eqn:laplace_bem_final}
	\frac{1}{2}\phi + \int_{\Gamma} \phi\partiald{G}{\nhat}\;\di{\Gamma} = \int_{\Gamma}\partiald{\phi}{\nhat}G\;\di{\Gamma},
\end{equation}

\noindent where $G = 1/4\pi r$ is the free-space Green's function for the Laplace equation ($\nabla^{2}G = -\delta$),  $\partiald{\cdots}{\nhat}$ represents the partial derivative in the direction normal to the boundary surface, and the (Cauchy principal-value) integrals are on the boundary $\Gamma$ of the domain. (The details of the derivation can be found in Ref. \cite{BrebbiaDominguez1992}.) The boundary element method (\bem) consists of discretizing the boundary into surface panels and enforcing Equation \eqref{eqn:laplace_bem_final} on a set of target points (collocation version). In a typical \bem, surface panels take a constant value $\phi_j$, and the surface integrals become sums over $N$ flat surface elements, $\Gamma_j$, resulting in the following discretized equation:
\begin{equation}
	\frac{1}{2}\phi_i = \sum_{j, i \neq j}^{N} \partiald{\phi_j}{\nhat_j}\;\int_{\Gamma}G_{ij}\di{\Gamma_j} - \sum_{j, i \neq j}^{N} \phi_j\int_{\Gamma}\partiald{G_{ij}}{\nhat_j}\;\di{\Gamma_j}.
\end{equation}

The values of the potential or its normal derivative on each panel are known from boundary conditions, resulting in either first-kind or second-kind integral equations. Finding the remaining unknowns requires solving a system of linear equations $A\vect{x}=\vect{b}$, where the elements of the coefficient matrix are
\begin{equation} \label{eqn:laplace_matrix}
	A_{ij} = 
	\begin{cases}
		\int_{\Gamma} G_{ij}\;\di{\Gamma_j}, & \phi\;\text{given on panel}\;j \\
		\int_{\Gamma} \partiald{G_{ij}}{\nhat_j}\;\di{\Gamma_j}, & \partiald{\phi}{\nhat}\;\text{given on panel } j
	\end{cases}
\end{equation}

\noindent
and $\vect{b}$ is formed with the known terms: e.g., if $\phi$ is given on panel $j$, then $\phi_j\int_{\Gamma_j}\partialdi{G_{ij}}{\nhat_j}\;\di{\Gamma_j}$ will appear in the term $b_i$ on the right-hand side of the linear system.

Expressing $G_{ij}$ and $\partialdi{G_{ij}}{\nhat_j}$ in terms of $1/r$ and $\nhat_j\cdot\nabla(1/r)$
\begin{eqnarray}
	\label{eqn:laplace_bem_G}\int_{\Gamma} G_{ij}\;\di{\Gamma_j} & = & \int_{\Gamma} \frac{1}{|\vect{x}_i-\vect{x}_j|} \;\di{\Gamma_j} \\ 
	\label{eqn:laplace_bem_dGdn}\int_{\Gamma} \partiald{G_{ij}}{\nhat_j}\;\di{\Gamma_j} & = & \int_{\Gamma}\frac{d\vect{x}\cdot\nhat_j}{|\vect{x}_i-\vect{x}_j|^{3}}\;\di{\Gamma_j}
\end{eqnarray}

The next steps are to apply an appropriate numerical integration scheme in order to generate all the terms of the coefficient matrix, and subsequently solve the linear system of equations, as described below.

\subsection{Boundary-integral solution of the Stokes equation}

The Stokes equation for a flow at very low Reynolds number, $\mu\nabla^{2}\vect{u} =  \nabla p$ (where $\mu$ is the viscosity and $p$ the pressure), can be rewritten in its integral formulation by means of a similar process as that described above for the Laplace equation. But it is a vector equation and its fundamental solutions are tensors. The boundary integral form of the the Stokes equation is

\begin{equation}
	\label{eqn:stokes_bem_12}
	\frac{1}{2}u_j(\vect{x_0}) = -\frac{1}{8\pi\mu}\int_{\Gamma} t_i(\vect{x})G_{ij}(\vect{x},\vect{x}_0)\;\di{\Gamma} + \frac{1}{8\pi} \int_{\Gamma} u_i(\vect{x})T_{ijk}(\vect{x},\vect{x}_0)n_k(\vect{x})\;\di{\Gamma}.
\end{equation}

\noindent where $\vect{u}$ is the velocity vector satisfying the Stokes equation (Einstein indicial summation implied), with $\sigma$ the corresponding stress tensor and $\vect{t} = \sigma\cdot\nhat$  the traction, vectors $\vect{x}$ and $\vect{x}_0$ are two distinct points in the domain, and $\vect{G}$ and $\vect{T}$ are the stokeslet and stresslet fundamental solutions:
\begin{equation}
	\label{eqn:stokeslet}
	G_{ij}(\vect{x},\vect{y})  =  \frac{\delta_{ij}}{r} + \frac{(x_i-y_i)(x_j-y_j)}{r^{3}} 
\end{equation}

\begin{equation}
	\label{eqn:stresslet}
	T_{ijk}(\vect{x},\vect{y},\nhat)  =  6\frac{(x_i-y_i)(x_j-y_j)(x_k-y_k)n_k}{r^{5}} \end{equation}

Indices $i, j, k$ denote here the Cartesian-tensor components and $\delta_{ij}$ is the Kronecker delta. Discretizing the boundary with $N$ surface panels results in sums that we now number with the index $J$.
The discretized form with constant surface panels becomes
\begin{equation}
	\label{eqn:stokes_bem_discretized}
	\frac{1}{2}u_j(\vect{x_0}) = -\frac{1}{8\pi\mu}\sum_{J=1}^{N}t_i\int_{\Gamma} G_{ij}(\vect{x}_J, \vect{x}_0)\;\di{\Gamma_J} + \frac{1}{8\pi} \sum_{J=0}^{N}u_i\int_{\Gamma} T_{ijk}(\vect{x}_J, \vect{x}_0)\cdot n_k\;\di{\Gamma_J}.
\end{equation}

\subsection{Numerical, semi-analytical and analytical integration methods}

\begin{figure}[t]
	\begin{centering}
\includegraphics[natwidth=5.15in,natheight=2.6in,width=0.45\textwidth]{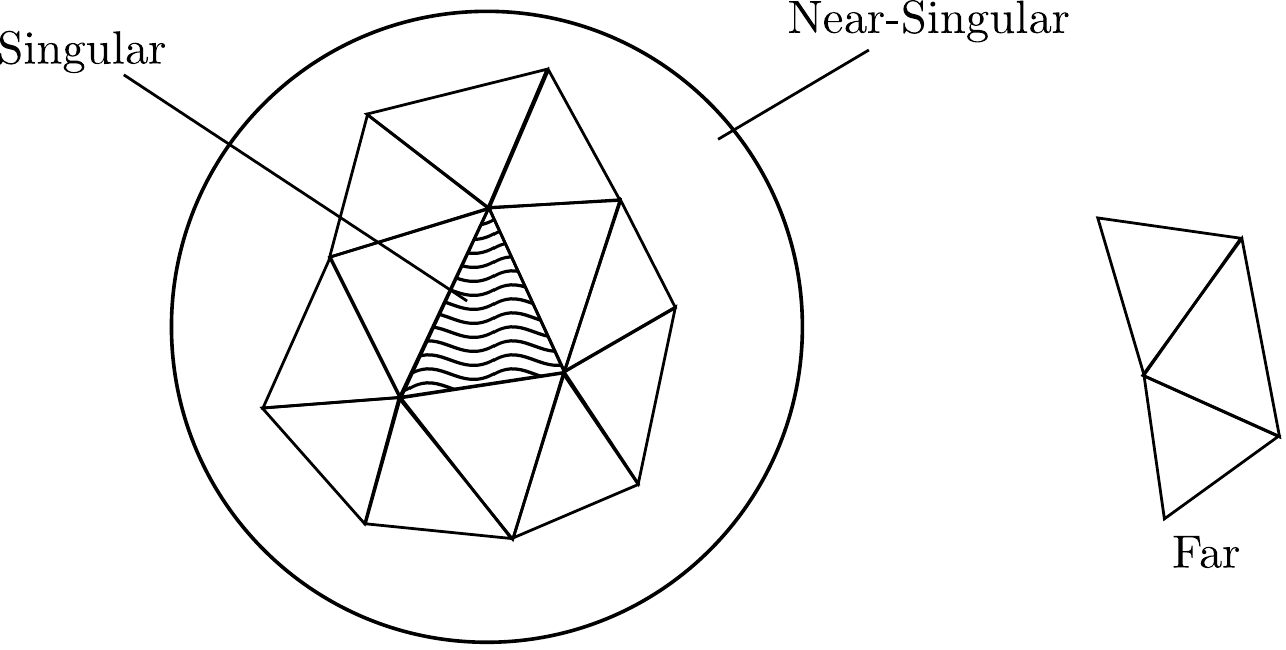}
	\caption{Depending on the distance between a source panel and target point, we use different numerical, analytical or semi-analytical integration methods, balancing computational efficiency and accuracy near singularities. The threshold distance between far and near-singular regions is $2\sqrt{2 S_j}$, where $S_j$ is the surface area of the source panel.}
	\label{fig:integration_domain}
	\end{centering}
\end{figure}

The boundary integral formulations all demand that we compute integrals of the type $\int_{\Gamma} \K_{ij}\;\di{\Gamma_j}$, where $\K_{ij}=\K(\vect{x}_j-\vect{x}_i)$ is the kernel, the point $\vect{x}_j$ is on the panel surface $\Gamma_j$ and the point $\vect{x}_i$ is a target or evaluation point. Because the kernel $\K$ is often singular, we need specific approaches depending on the distance $\vect{x}_j-\vect{x}_i$. Where the target point is far enough from the surface $\Gamma_j$, simple quadrature with a few Gauss points will suffice. As the target point \emph{nears} the source panel (with the definition of ``near'' to be determined), we need high-accuracy quadrature. Finally, in the case where the target point is on the source panel, the integral is (close to) singular and we must use analytical or semi-analytical methods. Figure \ref{fig:integration_domain} illustrates the three situations.

Applying Gauss quadrature to the integrals appearing in the coefficients of the boundary element discretization of the Laplace equation, \eqref{eqn:laplace_bem_G} and \eqref{eqn:laplace_bem_dGdn}, for example, gives
\begin{eqnarray}
	\label{eqn:gauss:1st-kind}
	\int_{\Gamma} G(\vect{x}_i,\vect{x}_j)\;\di{\Gamma_j} & \approx & \sum_{k=1}^{K} q_k\cdot S_j\cdot \frac{1}{|\vect{x}_i-\vect{x}_k|},\;\;\vect{x}_k \in \Gamma_j, \\ 
	\label{eqn:gauss:2nd-kind}
	\int_{\Gamma} \partiald{G(\vect{x}_i,\vect{x}_j)}{\nhat_j}\;\di{\Gamma_j} & \approx & \sum_{k=1}^{K}q_k\cdot S_j\cdot \frac{d\vect{x}\cdot\nhat_j}{|\vect{x}_i-\vect{x}_k|^{3}},\;\;\vect{x}_k \in \Gamma_j,
\end{eqnarray}

\noindent
with $q_k$ the area-normalized Gauss quadrature weights and $S_j$ the surface area of panel $\Gamma_j$. To control the accuracy of the numerical integration, we vary the number of quadrature points, using for example $K= 4$ for targets that are far from the source panel and $K\approx 20$ for near-singular situations. 
The near-singular region is within a distance of $2\sqrt{2 S_j}$, a criterion that we settled on after  testing with several choices using a panel's characteristic length scale as factor.
When the target point is on the source panel, the standard approach is to use analytical or semi-analytical methods for the singular and hyper-singular integrals over the panels.
For the Laplace equation, we used a semi-analytical method. 
It applies the technique first presented in the classic work of Hess and Smith \cite[p.~49, ff.]{HessSmith1967} for decomposing the integral into a sum of integrals over triangles formed by the projection of the target point on the panel plane, and the panel vertices. Using polar coordinates, the integration over the radial component can be done analytically and the integration over the angular component is done by quadrature. 
Several analytic integration techniques are at our disposal for dealing with the singular integrals from boundary element methods. Explicit expressions for these integrals over flat triangular domains result in recursive formulae on the edges of the integration triangle. These are available for Laplace potentials \cite{Fata2009} and linear elastic surface potentials \cite{Fata2011}. 
We obtained the analytic integrals for the Stokes equation from Fata's formulas for linear elasticity, after setting the Poisson ratio to $1/2$. 

\subsection{Krylov subspace methods}

For large linear systems of equations $A\vect{x}=\vect{b}$, direct solution is generally impractical and iterative solution methods are preferred. Krylov subspace methods derive from the Cayley-Hamilton theorem, which states that you can express the inverse of a  matrix $A$ as a linear combination of its powers. A Krylov subspace is spanned by the products of $b$ and powers of $A$; to order-$r$, this is: $K_{r}(A,b) = \text{span}\{ b, Ab, A^{2}b, ..., A^{r-1}b\}$.
Krylov methods include the conjugate gradient method, the biconjugate gradient stabilized method (\textsc{bicgstab}) and the generalized minimal residual method, \gmres \cite{SaadSchultz1986}, which we use. The greatest cost per iteration is the matrix-vector product (mat-vec), $w\gets A\cdot z_j$, taking $\O{N^{2}}$ time in a direct implementation. However, given the structure of the coefficient matrix in boundary element methods, this operation can be reduced to $\O{N}$ time using, for example, the fast multipole method.
The key is understanding the matrix-vector product as an $N$-body problem. Let's consider the first-kind integral problem for the Laplace equation. In \eqref{eqn:laplace_matrix} we see that the matrix coefficients are $\int_{\Gamma} G_{ij}\,\di{\Gamma_j}$. Applying Gauss quadrature to obtain the coefficients gives \eqref{eqn:gauss:1st-kind}. Thus, the matrix-vector product gives the following  for row $i$
\begin{equation}\label{eqn:matvec-onerow}
	\sum_{j=1}^{N_p}  \partiald{\phi_j}{\nhat_j}\; \sum_{k=1}^{K} q_k\cdot S_j\cdot \frac{1}{|\vect{x}_i-\vect{x}_k|},\;\;\vect{x}_k \in \Gamma_j, 
\end{equation}

\noindent
The inner sum is over the Gauss quadrature points (only a few per panel) and the outer sum is over the integration panels. We list the pseudocode of the full mat-vec in the Appendix as Algorithm \ref{alg:matvec}. Taking all the quadrature points together as the set of ``sources'' and all the collocation points on the panels as the set of ``targets,'' the algorithm is reduced to two for-loops instead of three, making the analogy with an $N$-body problem more clear.

\subsection{Fast multipole method}

Fast multipole methods were invented to accelerate the solution of $N$-body problems, that is, problems seeking to determine the motion of $N$ bodies that interact with each other via a long-distance effect (like electrostatics or gravitation). A direct approach to such a problem takes $\O{N^{2}}$ time to compute. The first \emph{fast} algorithms for $N$-body problems \cite{Appel1985,BarnesHut1986} combined two ideas: (1) approximating the effect of groups of distant bodies with a few moments (of the charges or masses), and (2) using a hierarchical sub-division of space to determine the acceptable distances to apply these approximations.
 These ideas produced the treecode algorithm, with $\bigO(N\log N)$ time to compute.
The fast multipole method \cite{GreengardRokhlin1987} introduces a third key idea that leads to $\bigO(N)$ scaling: allowing groups of distant bodies to interact with \emph{groups} of targets, by means of a mathematical representation called local expansion.

\begin{figure}
\begin{center}
	\includegraphics[width=0.5\textwidth]{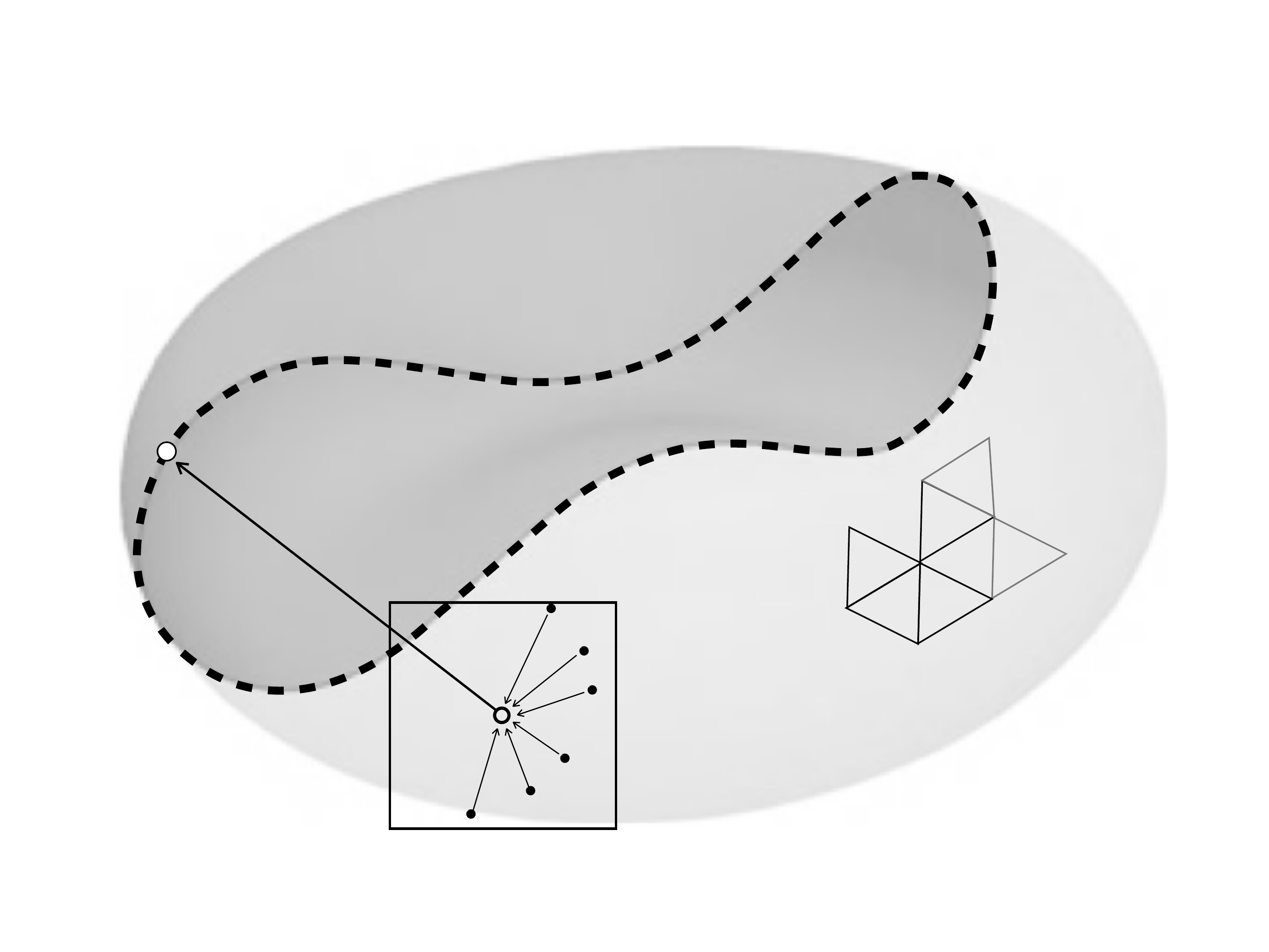}
	\caption{First step in viewing the influence of boundary elements as an $N$-body problem accelerated by the fast multipole method: Gauss quadrature points in a region of the surface triangulation are grouped together in a cluster, and their contribution to the potential at a remote target point (white circle) is considered via a series expansion at the cluster center.}
	\label{fig:rbc_fmmbox}
\end{center}
\end{figure}

A typical $N$-body problem evaluates a potential $\phi$ on $i=1, 2, \cdots, N$ bodies
using the following expression
\begin{equation}\label{eqn:nbody}
	\phi_{i} = \sum_{j=0}^{N} m_{j}\cdot\K(\vect{x}_{i},\vect{y}_{j}) \; = \; \sum_{j=0}^{N}\K_{ij}m_{j},
\end{equation}

\noindent where $\K_{ij} = \K(\vect{x}_{i},\vect{y}_{j})$ is referred to as the \emph{kernel}, and the potential is a solution of an elliptic equation, e.g., the Poisson equation $\vect{F}_i = - \nabla^2 \phi_i$ for gravitation. The expression in \eqref{eqn:nbody} is analogous to that for one row of the \bem\ mat-vec in Equation \eqref{eqn:matvec-onerow}, taking all $N_p \cdot K$ Gauss points collectively.
The first step in the \fmm\ acceleration of \bem\ is to group the Gauss quadrature points (i.e., the ``sources'') into clusters, and represent their influence via multipole expansions at the cluster centers. If using Taylor expansions, for example, truncated to the first $p$ terms, the potential at a point is approximated by
\begin{equation}
	\phi(\vect{x}_i) \approx \sum_{||\vect{k}||=0}^{p}\frac{1}{\vect{k}!}D^{\vect{k}}_{\vect{y}} \K(\vect{x}_i,\vect{y}_c)\, \sum_{j=1}^{n_c} m_j (\vect{y}_j-\vect{y}_c)^{\vect{k}},
	\label{eqn:cartesian_multipole}
\end{equation}

\noindent where $\vect{k}=\{k_1, k_2, k_3\}$ is a multi-index, $\vect{k}! = k_1!k_2!k_3!$, $\vect{y}^{\vect{k}} = y_1^{k_1}y_2^{k_2}y_3^{k_3}$, $D_{\vect{y}}^{\vect{k}} = D^{k_1}_{y_1}D^{k_2}_{y_2}D^{k_3}_{y_3}$ is the derivative operator and $p$ is the order of the expansion. The right-most terms are the multipoles, i.e., powers of distances with respect to the cluster center, and they can all be pre-calculated for a set of sources. Forming the clusters consists of recursively sub-dividing the spatial domain until a limit number of sources $\ncrit$ remains in each box, resulting in an octree structure. 
The step of computing the multipole expansions for each source cluster at the deepest level of the tree is referred to as the particle-to-multipole operation, {\ptom}. These multipole expansions are then translated and added to represent larger clusters (at higher levels of the tree) in the multipole-to-multipole operation, {\mtom}. The process just described is called the \emph{upward sweep}.
At this point, a treecode algorithm evaluates the potential on target points, performing  multipole-to-particle operations, {\mtop}---in the \bem, the multipoles represent clusters of quadrature points and the potential is evaluated on collocation points, as illustrated in Figure \ref{fig:rbc_fmmbox}. As implied in Equation\eqref{eqn:cartesian_multipole}, this is an approximation; and as implied in Figure \ref{fig:rbc_fmmbox}, the approximation is acceptable for remote target points only. The parameter that dictates whether we make the approximation is the \emph{multipole acceptance criterion}, \mac, denoted by $\theta$ and enforced as the maximum allowed ratio between cluster size and distance between targets and cluster center. When the \mac\ is not satisfied, sources and targets interact directly via \eqref{eqn:nbody} (called particle-to-particle operation, \ptop). The key to achieving the optimal $\bigO(N)$ scaling are local expansions, representing a group of target points; thus, the \fmm\ adds three operations to the algorithm: the multipole-to-local transformation (\mtol), the local-to-local translation (\ltol) and the local-to-particle evaluation (\ltop).

The Cartesian-expansion \mtol\ operation scales as $\O{p^{6}}$, which quickly becomes expensive when requiring higher accuracy via higher $p$. Using spherical-harmonic expansions reduces this cost to  $\O{p^{4}}$, but the mathematical derivations and expressions become more cumbersome. For the sake of brevity, we omit all the details and write the final form for the spherical expansions of the Laplace potential---denoted here by $\Phi(\vect{x}_i)$ to differentiate with the spherical co-ordinate $\phi$---as follows

\begin{eqnarray}
	\Phi(\vect{x}_i) & = & \sum_{n=0}^{p}\sum_{m=-n}^{n}\frac{Y^{m}_n(\theta_i,\phi_i)}{r_i^{n+1}}\underbrace{\left \{ \sum_j^{N}q_j\rho^{n}_jY^{-m}_n(\alpha_i,\beta_i)\right \} }_{M^{m}_n} \, \text{and}\\
	\Phi(\vect{x}_i) & = & \sum_{n=0}^{p}\sum_{m=-n}^{n}r_i^{n}Y^{m}_n(\theta_i,\phi_i)\underbrace{\left \{ \sum_j^{N}q_j\frac{Y^{-m}_n(\alpha_i,\beta_i)}{\rho^{n+1}_j}\right \} }_{L^{m}_n}.
\end{eqnarray}

\noindent
Here, $M^{m}_n$ and $L^{m}_n$ denote the multipole and local expansion coefficients; $Y_{n}^{m}$ is the spherical harmonic function; $(r,\theta,\phi)$ and $(\rho,\alpha,\beta)$ are the distance vectors from the center of the expansion to points $\vect{x}_i$ and $\vect{x}_j$, and $q_j$ are the weights of sources. For our purposes, we just want to make clear the parallel between a fast $N$-body algorithm and a fast \bem\ matrix-vector product via this necessarily abbreviated summary of the \fmm.

\subsection{Inexact Krylov iterations and relaxation strategies}

Accelerating a \bem\ solution with \fmm\ implies computing the matrix-vector products with some error, and the parameter controlling this accuracy is the order of multipole expansions, $p$. 
Two natural questions to ask are how large does the value of $p$ need to be to ensure that the iterative solver will still converge and what's the impact on the accuracy of the converged solution due to the error on the computed mat-vec.
We might expect to choose that value of $p$ that ensures convergence and desired accuracy, and apply it evenly for all iterations. 
But it turns out that Krylov iterations have a surprising property: the \emph{first} iteration needs to be computed with high accuracy, but accuracy requirements can be \emph{relaxed} for later iterations. 
Translation operators in the \fmm can scale as $\O{p^{6}}$ or $\O{p^{4}}$---using Cartesian or spherical expansions, respectively---, so this property of Krylov methods could offer the potential for further accelerating the \bem\ solution.
Bouras and Frayss{\'e} \cite{bouras2000relaxation,bourasfraysse2005} studied the effect of inexact Krylov iterations on convergence and accuracy via numerical experiments.
They found that if the system matrix is perturbed, so we are computing $(A+\Delta A_k)\vect{z}$ on each iteration, and the perturbation stays nearly equal in norm to $\eta \|A\|$, then the computed solution will have an error of the same order, $\eta$. This is the situation when simply computing the mat-vec within the limits of machine precision. But they also showed the more surprising result that the magnitude of the perturbations $\Delta A_k$ can be allowed to grow as the iterations progress.
They define a relaxation strategy whereby, for a desired final tolerance $\eta$ in the solution of the system, the mat-vecs in each iteration are computed with a coefficient matrix perturbed with $\Delta A_k$, where if $r_k$ is the residual at step $k$,
\begin{equation}\label{eqn:matrix-perturbation}
  \|\Delta A_k\| = \varepsilon_k\|A\|, \quad \varepsilon_k=\min\left( \frac{\eta}{\min(\|r_{k-1}\|,1)}, 1\right).
\end{equation}

\noindent The matrix perturbation is always larger than the target tolerance $\eta$, and $\varepsilon_k$ increases when $r_k$ decreases (without surpassing $1$). In other words, the accuracy of the system mat-vecs are \emph{relaxed} as iterations proceed. This relaxation strategy led to converged solutions with \gmres, \textsc{cg} and \textsc{bicgstab}, using a variety of test matrices---often, and remarkably, in about the same number of iterations as a non-relaxed solution. In sum, Krylov methods proved to be robust to inexact mat-vecs and only the first Krylov vectors need to be computed accurately.
The numerical evidence is also supported by theoretical studies \cite{simonciniszyld2003,vandeneshofsleijpen2004}.

In the fast multipole method, we have error bounds available for the approximations made in the various operations that make up the algorithm. Using the spherical-harmonics expansion for the Laplace kernel, for example, the error is bounded as follows 
\begin{equation}\label{eqn:multipole_error}
	\left | \phi(r, \theta, \phi) - \sum_{n=0}^{p}\sum_{m=-n}^{n}\frac{M^{m}_{n}}{r^{n+1}}\cdot Y^{m}_{n}(\theta, \phi) \right | \leq \frac{\sum_{i=1}^{N}q_{i}}{r-a}\left ( \frac{a}{r} \right )^{p+1},
\end{equation}

\noindent where $a$ isthe  cluster radius and $r$ is the distance between the multipole center and the target. The above inequality is given in Ref. \cite[p.~55]{greengard1987} with label (3.38), and proved using the triangle inequality. This reference also gives similar bounds for the translation and evaluation operators. In the particular case that $r/a\geq 2$ (the distance between a multipole center and target is at least twice the cluster radius), the number of terms needed to obtain a given accuracy $\varepsilon$ is 

\begin{equation}\label{eqn:fmm_p}
	p \sim \lceil -\log_{2}(\varepsilon) \rceil.
\end{equation}

\noindent In combination with Equation \eqref{eqn:matrix-perturbation}, we thus have an explicit  relation between the allowed perturbation on the mat-vec for the inexact Krylov iterations to converge, and the order of the multipole expansion in the \fmm.
Although the \fmm error bounds are known to be rather loose, we use this approach for a conservative relaxation of $p$, and we expect that with experience and more detailed studies, the relaxation strategy could be pushed to give better speed-ups in specific applications.

\section{Results and discussion}
\begin{figure}[h]
\begin{center}
	\includegraphics[width=0.5 \textwidth]{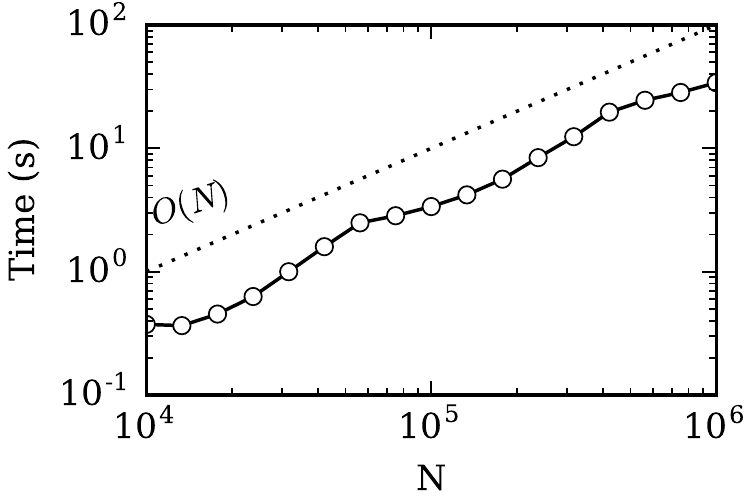}
	\caption{Scaling of the \fmm\ code with respect to problem size $N$, using a Laplace kernel with spherical expansions, $p=5$ and  $\ncrit = 126$. The tests ran on a lab workstation using six \cpu\ threads. The dotted line shows the expected $\O{N}$ scaling. Plotting script and figure available under CC-BY \cite{WangLaytonBarba2016-figshare1}.}
	\label{fig:fmm_scaling}
\end{center}
\end{figure}

Our results include numerical experiments with \bem\ solutions for the Laplace and Stokes equations, including an application to Stokes flow around red blood cells. The \bem\ solver is accelerated with a fast multipole method using spherical expansions, implemented in a code designed to use boundary elements (panels) as the sources and targets, and written in multi-threaded \cpp. To facilitate reproducibility of our research, the code is made available via its version-control repository,\footnote{\href{https://github.com/barbagroup/fmm-bem-relaxed}{https://github.com/barbagroup/fmm-bem-relaxed}} under an MIT license. The paper manuscript is also available via its repository,\footnote{\href{https://github.com/barbagroup/inexact-gmres}{https://github.com/barbagroup/inexact-gmres}} which includes running scripts and plotting scripts to produce the included figures (see each figure caption for details).
We ran all our tests on a lab workstation with an Intel Core i7-5930K hexa-core processor and 32GB RAM. Using a test with the Laplace kernel, we confirmed that our \fmm\ code scales as $\O{N}$ by timing several runs with increasing problem size; see Figure \ref{fig:fmm_scaling}. We also found that the preconditioned cases take fewer iterations but more time to converge than the unpreconditioned ones, because the residuals after the first several iterations in the preconditioned cases are greater. For an inexact {\gmres} iteration, a bigger residual leads to a higher required $p$, which offsets the benefit of using fewer iterations. Therefore, we performed these tests without preconditioning.
The subsequent experiments investigate the use of a Krylov relaxation strategy by reducing the order of multipole expansions, $p$, as the iterations progress to convergence, determining the speed-up provided by such a strategy for different scenarios and problem sizes.

\subsection{Inexact {\small GMRES} for the solution of the Laplace equation}
\label{sec:inexactLaplace}
\begin{figure}
\begin{center}
	\subfloat[][128 panels]{\includegraphics[natwidth=4.73in,natheight=3.94in,width=0.3\textwidth]{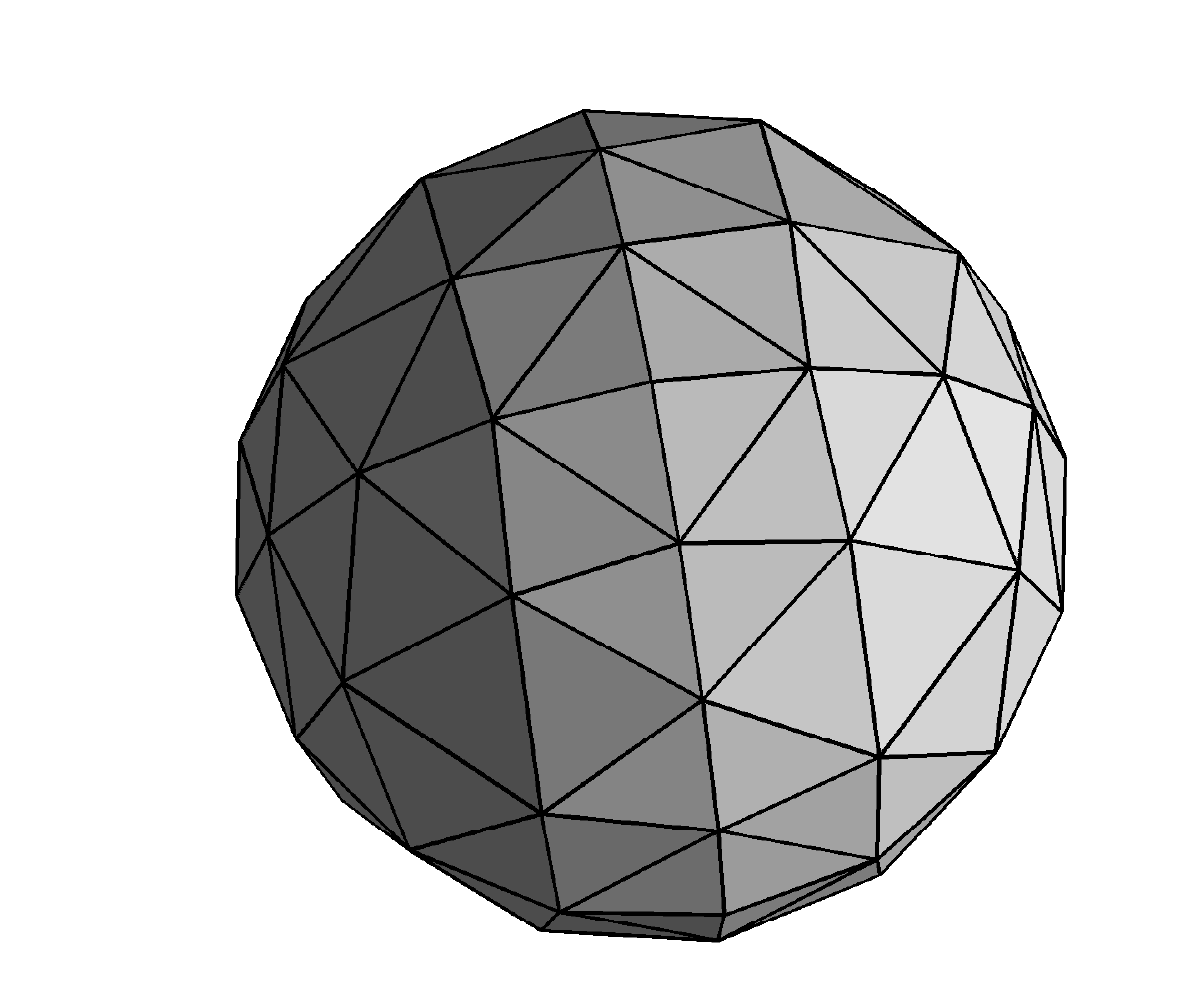}\label{fig:sphere128}}\qquad
	\subfloat[][2048 panels]{\includegraphics[natwidth=4.73in,natheight=3.94in,width=0.3\textwidth]{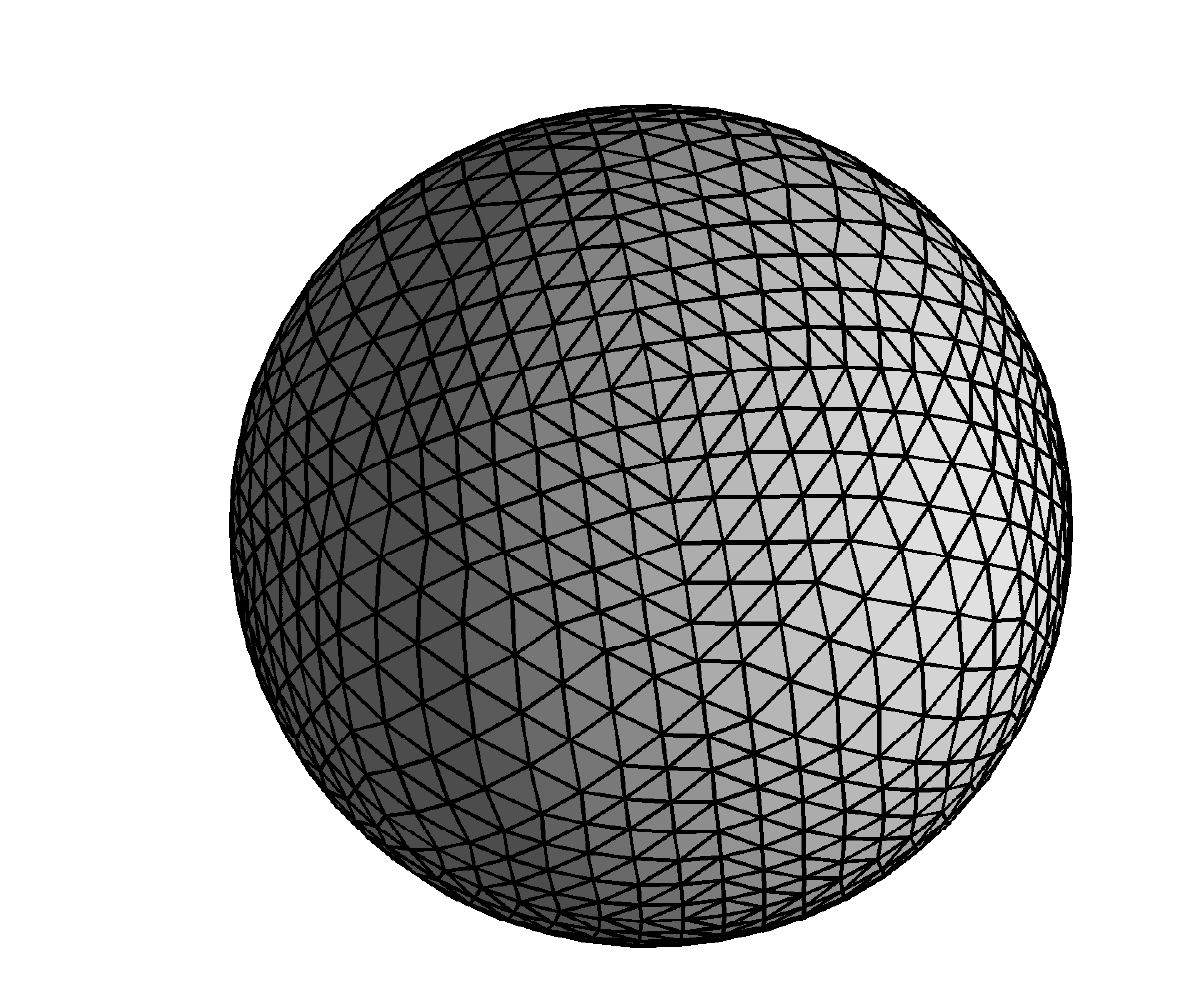}\label{fig:sphere2048}}
	\caption{Triangular discretizations of a spherical surface.}
	\label{fig:glob_spheres}
\end{center}
\end{figure}

\begin{figure}[t]
\begin{center}
	\includegraphics[natwidth=3in,natheight=2in,width=0.5\textwidth]{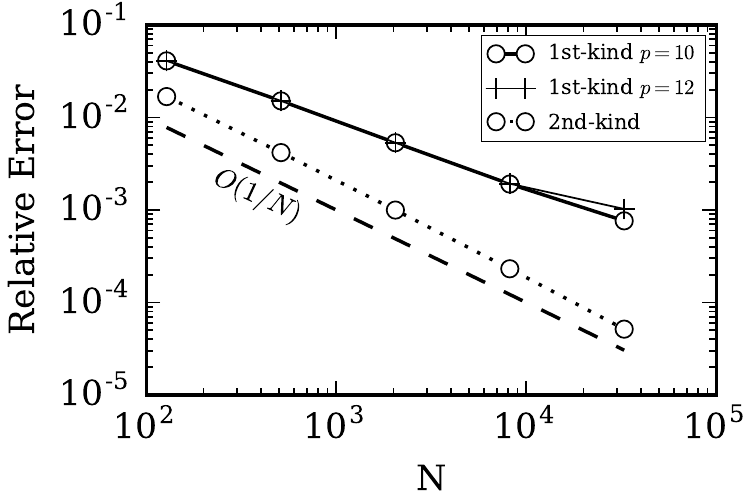}
	\caption{Convergence of 1st-kind (solid line) and 2nd-kind (dotted line) solvers for the Laplace equation on a sphere, using a \gmres with \fmm-accelerated matrix-vector products, with parameters: $\theta_{\text{MAC}} = 0.5$, $p=12$, $k=4$ and solver tolerance of $10^{-6}$ (no relaxation). For the 1st-kind formulation, we also used $p=10$ (thinner solid line) to show a degradation of grid-convergence with lower $p$. The vertical axis is the $L^2$-norm of the relative error with respect to the analytical solution for a constant potential or charge on the surface: $\phi = \partialdi{\phi}{\nhat} = 1$. Plotting script and figure available under CC-BY \cite{WangLaytonBarba2016-figshare2}.}
	\label{fig:laplaceconvergence}
\end{center}
\end{figure}

To start, we studied grid convergence comparing numerical results with the analytical solution using a sphere with constant potential and charge on the surface: $\phi = \partialdi{\phi}{\nhat} = 1$. To make surface triangulations of a sphere with increasing refinement, we started with an 8-triangle closed surface, then split recursively each triangle into four smaller ones. Figure \ref{fig:glob_spheres} shows two example discretizations. We solved the boundary-element problem by collocation in both the first-kind and second-kind integral formulations, using a standard \gmres with fast-multipole-accelerated mat-vecs and the semi-analytical integrals for the singular terms. For the far-field approximations, we used spherical-harmonic expansions with the following parameters in the \fmm: $\theta_{\text{MAC}} = 0.5$, $p = 12$; the number of Gauss points was $k=4$ and the tolerance in the iterative solver was $10^{-6}$. 
Figure \ref{fig:laplaceconvergence} shows the resulting convergence for both first-kind and second-kind formulations of the boundary element method on a sphere. The observed order of convergence is 0.76 for the 1st-kind formulation and 1.02 for the 2nd-kind formulation, computed with the three points in the middle of each line. We also observed that the finest mesh detracts from the convergence line with $p=10$, as the thinner line shows in Figure \ref{fig:laplaceconvergence}; in that case, the discretization error of a simple geometry like a sphere is very small and the error of the \fmm-accelerated mat-vecs overtakes.
This convergence analysis gives confidence in our \bem code, the singular/near-singular integral calculations, and the far-field approximation using the \fmm.


Next, we looked at the following test to see how the residual changes as the \gmres iterations proceed and  what value of $p$ is required in the \fmm-accelerated mat-vecs to continue convergence, according to Equation \eqref{eqn:fmm_p}. We discretized a sphere with $32,768$ surface triangles and solved a first-kind integral equation using a solver tolerance of $10^{-6}$ with an initial value of $p$ set to 8 and 12. As the residual gets smaller, the value of $p$ needed to maintain convergence of the solver drops, and a low-$p$ of just 2 is sufficient by the ninth iteration for $p_\text{{initial}}=8$ and by the fourth iteration for $p_{\text{initial}}=12$. This offers the potential for substantial speed-ups in the calculations, because the translation operators of the \fmm scale from $\bigO(p^{4})$ for spherical harmonics to $\bigO(p^{6})$ for Cartesian expansions. A more accurate first iteration results in a faster drop of required-$p$ and fewer iterations. 
But we note that only the far-field evaluation can be sped-up with the relaxation strategy, which means that the correct balance between near field and far field in the \fmm could change as we reset $p$ in the later iterations.

\begin{figure}
	\centering
	\includegraphics[natwidth=7in,natheight=2in,width=0.95\textwidth]{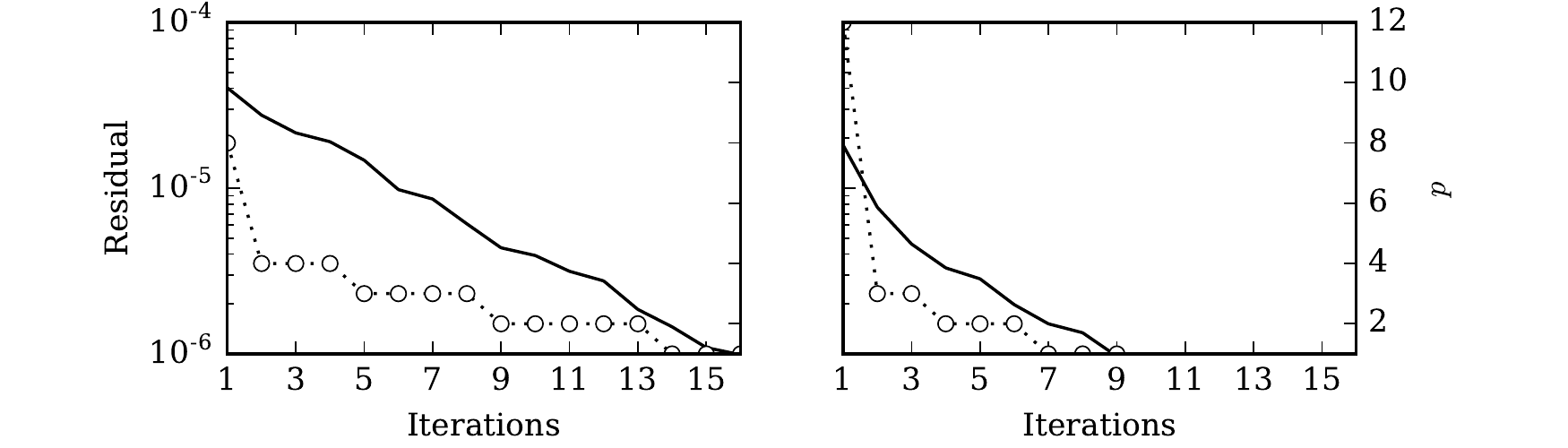}
	\caption{In a test using a sphere discretized with $32,768$ triangles, with $k=4$, $p_\text{initial}=8$ (left plot) and $p_\text{initial}=12$ (right plot), the residual $\|r_{k}\|$  (solid line, left axis) decreases with successive \gmres iterations while the necessary $p$ (open circles, right axis) to achieve convergence drops quickly. Plotting script and figure available under CC-BY \cite{WangLaytonBarba2016-figshare2}.}
	\label{fig:residualp}
\end{figure}

\begin{figure}
	\centering
	\includegraphics[natwidth=3in,natheight=2in,width=0.45\textwidth]{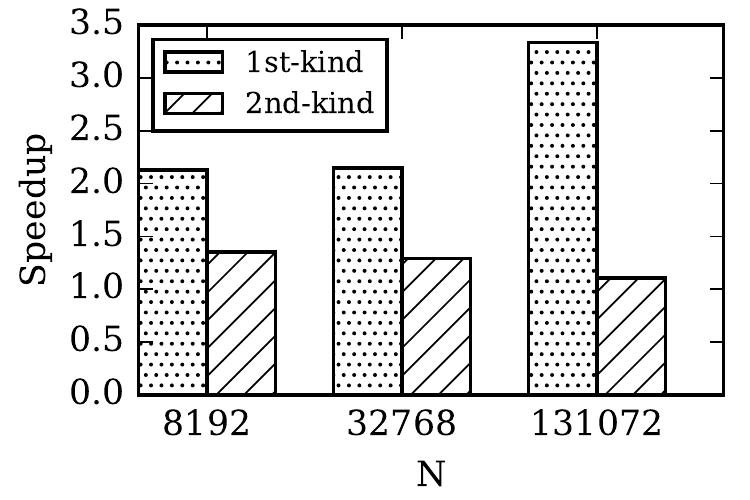}
	\caption{Speed-up using a relaxation strategy for three different triangulations of a sphere ($N$ is the number of surface panels), using 1st-kind and 2nd-kind integral formulations. Initial $p=10$, $k=4$, solver tolerance $10^{-6}$. Time is measured by averaging the solving time of three identical runs. (Multi-threaded evaluator running on 6 \cpu\ cores.) Plotting script and figure available under CC-BY \cite{WangLaytonBarba2016-figshare2}.}
	\label{fig:relaxation_timing}
\end{figure}

To find out the potential speed-up, we compared the time to solution for different cases with and without the relaxation strategy. Using three surface discretizations, we solved the boundary-element problem with 1st- and 2nd-kind formulations, with a multi-threaded evaluator on 6 \cpu\ cores. In each case, we were careful to set the value of $\ncrit$ to minimize the time to solution of the particular test case. 
The detailed results are given in Tables \ref{tab:laplace_1st_relaxation} and \ref{tab:laplace_2nd_relaxation}.
Figure \ref{fig:relaxation_timing} shows the speed-up in the time spent solving the linear system of equations to the specified tolerance. 
As indicated in the caption, we used an initial value of $p=10$ and a solver tolerance of $10^{-6}$.
We also repeated the experiments with less stringent accuracy settings of $p=8$, and solver tolerance $10^{-5}$. In that case, the speedups are reduced (maximum observed speedup of 2.22 for the 1st-kind formulation), but we decided to show the higher accuracy tests, given the observed degradation of grid-convergence with lower $p$, as seen in Figure \ref{fig:laplaceconvergence}.

\begin{table}[h]
\footnotesize
\begin{center}
\begin{tabular}{c|cc|cc|cc}
  & \multicolumn{2}{c|}{Non-Relaxed} & \multicolumn{2}{c|}{Relaxed} &  & Number\\
  N & $\ncrit$ & $\tsolve$ & $\ncrit$ & $\tsolve$ & Speed-up & of iterations\\
 \hline
   & & & & & &\\
  8192 & 400 & 4.93 & 175 & 2.32 & 2.13 & 11\\
  32768  & 400 & 20.32 & 200 & 9.45 & 2.15 & 13\\
  131072  & 400 & 137.84 & 200 & 41.32 & 3.34 & 22\\
 
\end{tabular}
\end{center}
\caption{Speed-ups for the relaxation strategy on a Laplace 1st-kind integral solver, $p=10$, $k=4$, solver tolerance $10^{-6}$.}
\label{tab:laplace_1st_relaxation}
\end{table}%

\begin{table}[h]
\footnotesize
\begin{center}
\begin{tabular}{c|cc|cc|cc}
  & \multicolumn{2}{c|}{Non-Relaxed} & \multicolumn{2}{c|}{Relaxed} &  & Number \\
  N & $\ncrit$ & $\tsolve$ & $\ncrit$ & $\tsolve$ & Speed-up & of iterations \\
 \hline
   & & & & & & \\
  8192 & 400 & 2.10 & 300 & 1.55 & 1.35 & 3 \\
  32768 & 400 & 7.45 & 200 & 5.76 & 1.29 & 2 \\
  131072 & 400 & 23.63 & 200 & 21.32 & 1.11 & 2\\
 
\end{tabular}
\end{center}
\caption{Speed-ups for the relaxation strategy on a Laplace 2nd-kind integral solver, $p=10$, $k=4$, solver tolerance  $10^{-6}$.}
\label{tab:laplace_2nd_relaxation}
\end{table}%

The results on Tables \ref{tab:laplace_1st_relaxation} and \ref{tab:laplace_2nd_relaxation} show a speed-up of between $2.1\times$ and $3.3\times$ for the 1st-kind integral formulation and between $1.1\times$ and $1.35\times$ for the 2nd-kind formulation. 
The rightmost column shows the number of iterations to converge: we need 3 or less iterations for the 2nd-kind formulation, which explains why we have such low speedup compared with the 1st-kind formulation.
These tests also taught us that one has to give up on the idea of partitioning the domain between a near-field and a far-field in a way that balances the time spent computing each one---an accepted idea in \fmm applications. When relaxing the accuracy of the \gmres iterations, the time taken to compute the far field decreases significantly. This means that to minimize time-to-solution when using relaxed \gmres, the near and far fields should not be balanced, but rather the far field should be bloated. As a result, the first few iterations are completely dominated by the time to compute the far field, but this is offset by the benefit of much cheaper iterations from then on. This is a simple but unexpected and counter-intuitive algorithmic consequence of using inexact \gmres with \fmm.

To better evaluate the potential speedup in more general cases, we designed tests in the following two situations: (a) where higher accuracy of \fmm-accelerated mat-vec product is needed (necessitating an initially higher $p$) , and (b) where the linear system demands a greater number of iterations to reach a smaller desired tolerance. 

\begin{figure}
	\centering
	\includegraphics[natwidth=3in,natheight=2in,width=0.45\textwidth]{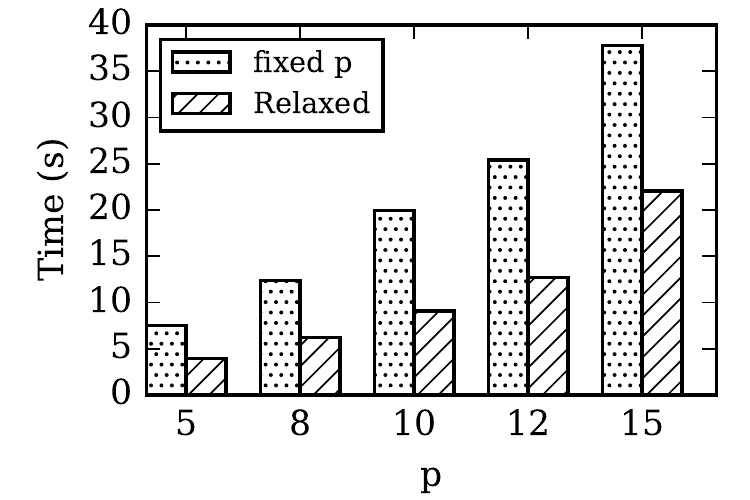}
	\caption{Timings for solving a 1st-kind Laplace integral formulation on a sphere discretized with $32,768$ panels, using a relaxed \gmres with different initial values of $p$, compared with a fixed-$p$ solver. The iteration count was capped at 10 for all cases. Time is measured by averaging the solving time of three identical runs. (Multi-threaded evaluator running on 6 \cpu\ cores.) Plotting script and figure available under CC-BY \cite{WangLaytonBarba2016-figshare2}.}
	\label{fig:laplace_p_speedup}
\end{figure}

Figure \ref{fig:laplace_p_speedup} shows the timings obtained in a situation where the initial $p$ is incrementally larger, representing applications that demand higher accuracy. To minimize the effect of discretization error on total accuracy, we discretized the sphere with $32,768$ panels. We also enforced a fixed number of $10$ \gmres iterations on this 1st-kind Laplace integral solver. Table \ref{tab:laplace_1st_p_relaxation} shows the data for this test: the speed-up is leveling off at around $2\times$ with a varying initial $p$. Note again that the shortest time-to-solution requires an unbalanced tree on the first iteration, with a bloated far field. This means that the first (high-$p$) iteration is much slower than the corresponding fixed-$p$ case: the speed-up is thus purely a product of the low-cost later iterations.

\begin{table}[h]
\footnotesize
\begin{center}
\begin{tabular}{c|cc|cc|c}
  & \multicolumn{2}{c|}{Relaxed} & \multicolumn{2}{c|}{Non-Relaxed} & \\
  $p$ & $\ncrit$ & $\tsolve$ & $\ncrit$ & $\tsolve$ & Speedup \\
   \hline
   & & & & & \\
  5 & 100 & 3.94 & 100 & 7.53 & 1.91 \\
  8 & 100 & 6.20 & 400 & 12.38 & 2.00 \\
  10 & 150 & 9.08 & 400 & 19.96 & 2.20 \\
  12 & 150 & 12.70  & 600 & 25.41 & 2.00 \\
  15 & 150 & 22.05 & 600 & 37.77 & 1.71 \\
 
\end{tabular}
\end{center}
\caption{Speed-up when using a relaxation strategy on a Laplace 1st-kind integral solver, compared with a non-relaxed solver, with increasing value of the initial $p$ (representing increased accuracy demands of the application), for a sphere discretized with $32,768$ panels. (Multi-threaded evaluator running on 6 \cpu\ cores.)}
\label{tab:laplace_1st_p_relaxation}
\end{table}%

\begin{figure}
	\centering
	\includegraphics[natwidth=3in,natheight=2in,width=0.5\textwidth]{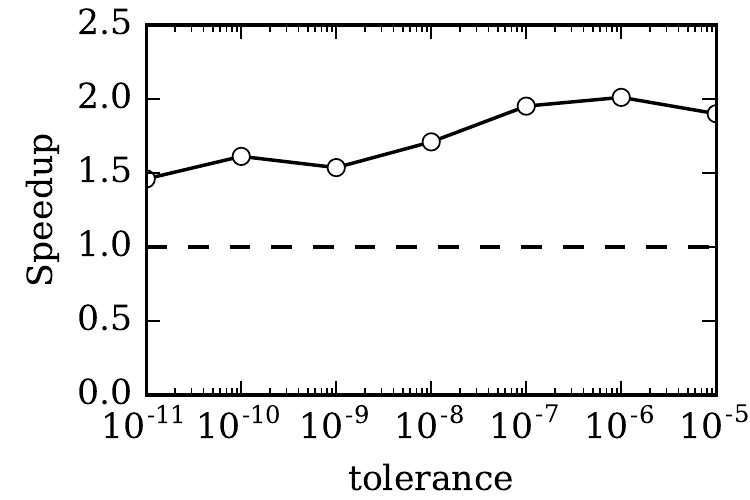}
	\caption{Speed-ups for solving a 1st-kind Laplace integral problem on a sphere discretized with $8,192$ panels, as the \gmres solver's tolerance increases; $p=10$ for all cases. Time is measured by averaging the solving time of three identical runs. (Multi-threaded evaluator running on 6 \cpu\ cores.) Plotting script and figure available under CC-BY \cite{WangLaytonBarba2016-figshare2}.}
	\label{fig:laplace_tolerance_speedup}
\end{figure}

The final case looks at the situation where the application might demand different tolerances. Keeping the value of initial $p$ fixed at 10, we ran several cases with $8,192$ panels and increasingly stricter tolerance. Figure \ref{fig:laplace_tolerance_speedup} and Table \ref{tab:laplace_1st_tolerance_relaxation} show that the relaxed cases experience a speed-up between $1.5\times$ and $2.0\times$ for different desired tolerances and the speed-up peaks at a tolerance of $10^{-6}$.

\begin{table}[h]
\footnotesize
\begin{center}
\begin{tabular}{c|cc|cc|c}
  & \multicolumn{2}{c|}{Relaxed} & \multicolumn{2}{c|}{Non-Relaxed} & \\
  tolerance & $\ncrit$ & $\tsolve$ & $\ncrit$ & $\tsolve$ & Speedup \\
 \hline
   & & & & & \\
  $10^{-5} $ & 400 & 1.65 & 400 & 3.13 & 1.90 \\
  $10^{-6} $ & 300 & 2.34 & 400 & 4.71 & 2.01 \\
  $10^{-7} $ & 300 & 4.07 & 400 & 7.93 & 1.95 \\
  $10^{-8} $ & 300 & 5.77 & 400 & 9.87 & 1.71 \\
  $10^{-9} $ & 300 & 8.14 & 400 & 12.51 & 1.54 \\
  $10^{-10} $ & 300 & 10.11 & 400 & 16.29 & 1.61 \\
  $10^{-11} $ & 300 & 12.24 & 400 & 17.88 & 1.46 \\
\end{tabular}
\end{center}
\caption{Speed-up of the relaxation strategy on solving a Laplace 1st-kind integral problem on a sphere discretized with $8,192$ panels, with a varying solver's tolerance; initial $p$ fixed to a value of 10. (Multi-threaded evaluator running on 6 \cpu\ cores.)}
\label{tab:laplace_1st_tolerance_relaxation}
\end{table}%

\subsection{Inexact {\small GMRES} for solving the Stokes equation}
Like in \S\ref{sec:inexactLaplace}, we start with a grid-convergence study to build confidence that the Stokes solver is correct and converges to the right solution at the expected rate. As an application of the Stokes equation, we chose low-Reynolds-number flow, using a spherical geometry for the grid-convergence study. This classical problem of fluid mechanics has an analytical solution that gives the drag force on the sphere as $F_d = 6\pi\mu Ru_x$, where $\mu$ is the viscosity of the fluid, $R$ is the Reynolds number and $u_x$ is the freestream velocity, taken in the $x$-direction. We solve a first-kind integral equation for the traction force, $\vect{t}$, by imposing $\vect{u} = (1,0,0)^{T}$ at the center of every panel, and compute the drag force with

\begin{equation}
	\label{eqn:stokes_traction_drag}
	F_d = \int_\Gamma t_x\;\di{\Gamma'} = \sum_{j=1}^{N} t_{x_j}\cdot A_j
\end{equation}

\begin{figure}
\begin{center}
	\includegraphics[natwidth=3in,natheight=2in,width=0.5\textwidth]{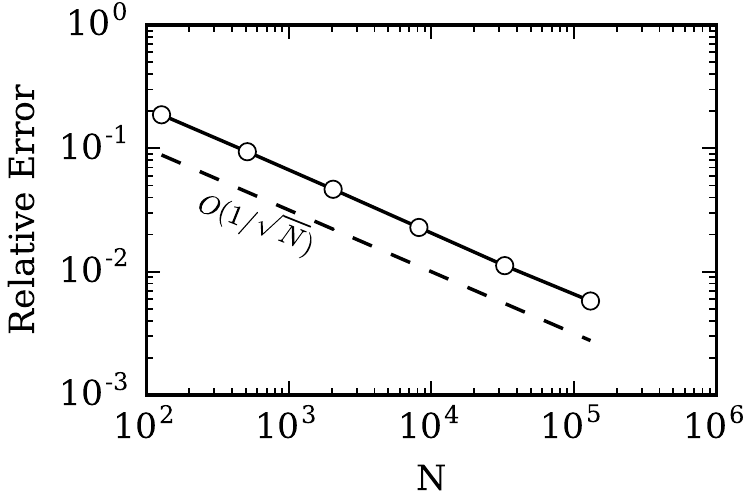}
	\caption{Convergence of the boundary-integral solution for Stokes flow around a sphere, using a  first-kind equation; $p=16$, linear system solved to $10^{-5}$ tolerance. The relative error is with respect to the analytical solution for drag on a sphere: $F_d = 6\pi\mu Ru_x$. Plotting script and figure available under CC-BY \cite{WangLaytonBarba2016-figshare3}.}
	\label{fig:stokes_convergence}
\end{center}
\end{figure}

For all the tests, we set $R=1$, $u_x = 1$ and $\mu = 10^{-3}$, giving a drag force of $F_d = 0.01885$. We solve the integral problem using a boundary element method with fast-multipole-accelerated mat-vecs in a \gmres solver. Based on the inexact Krylov method \cite{bouras2000relaxation,bourasfraysse2005}, we choose an initial $p$ of 16 to compute the first Krylov vectors with full accuracy, and use a $10^{-5}$ tolerance in the iterative solver.
Figure \ref{fig:stokes_convergence} shows that we observe convergence at the expected rate of $\O{1 / \sqrt{N}}$, for first-kind integral equations.

\begin{figure}
\begin{center}
	\includegraphics[natwidth=3in,natheight=2in,width=0.5\textwidth]{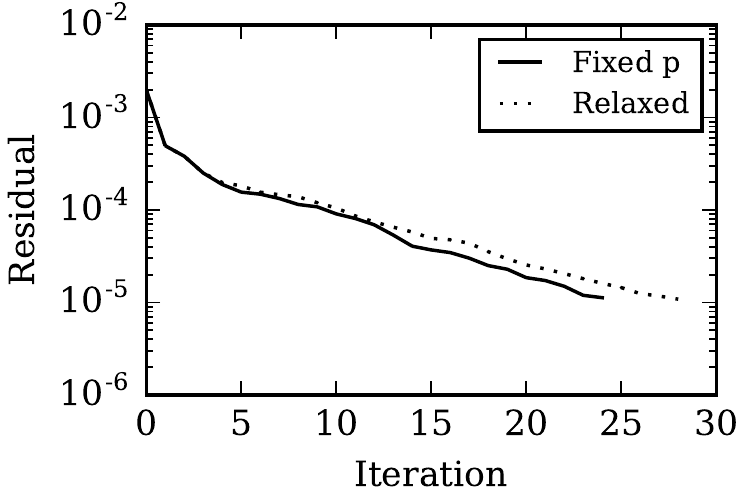}
	\caption{Residual history solving for surface traction on the surface of a sphere (first-kind integral problem), with a $10^{-5}$ solver tolerance, 8192 panels, and initial $p=16$ in the relaxed case (and throughout in the non-relaxed case). Plotting script and figure available under CC-BY \cite{WangLaytonBarba2016-figshare3}.}
	\label{fig:stokes_residual_history_relaxed}
\end{center}
\end{figure}

Like with the Laplace equation, we need to show that the relaxation strategy does not hinder convergence and that there is a potential for speed-ups. We solved the problem of Stokes flow around a sphere, discretized with $8,192$ panels, and compared the residual history of a fixed-$p$ solver with a relaxed \gmres with an initial $p=16$. Figure \ref{fig:stokes_residual_history_relaxed} shows that the residual history (for the traction force) is similar for both relaxed and non-relaxed Krylov iterations, both methods reaching the stipulated tolerance of $10^{-5}$ after about 27 iterations.
The number of iterations needed to converge is larger in the case of the Stokes equation compared to the Laplace equation, which bodes well for the speed-up that we could get from relaxation. Moreover, computing the spherical expansion for the Stokes kernel is equivalent to computing four Laplace expansions, which combines with the larger number of iterations to offer greater speed-ups.

\begin{figure}
\begin{center}
	\includegraphics[natwidth=9in,natheight=4in,width=0.95\textwidth]{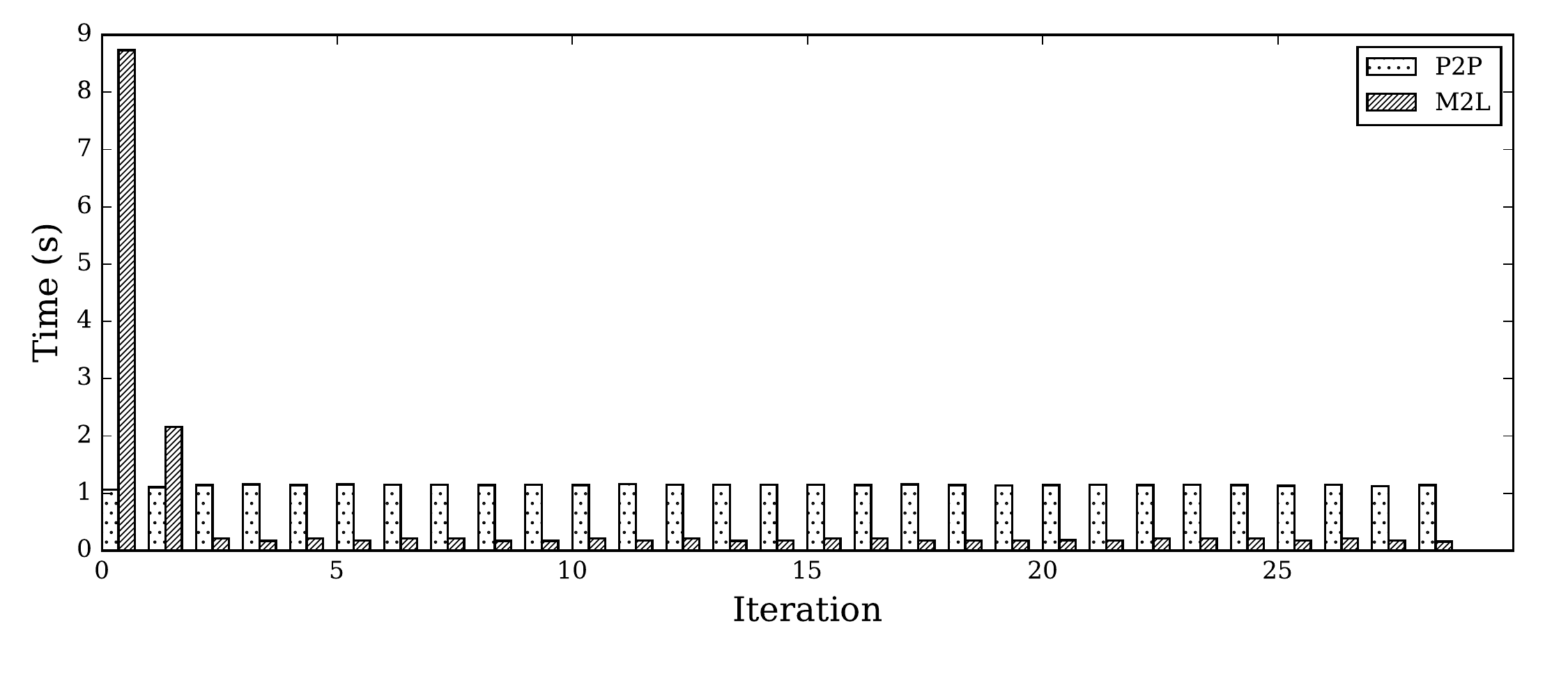}
	\caption{Time breakdown between {\ptop} and {\mtol} when using a relaxation strategy for solving surface traction on the surface of a sphere, $10^{-5}$ solver tolerance, $8,192$ panels, $p=16$. Plotting script and figure available under CC-BY \cite{WangLaytonBarba2016-figshare3}.}
	\label{fig:stokes_relaxation_breakdown}
\end{center}
\end{figure}

Figure \ref{fig:stokes_relaxation_breakdown} illustrates clearly how we need to adjust the balance between near field and far field when using relaxation strategies. Because most of the time is spent computing at the low values of $p$, we need to start with a bloated far field. The bar plot shows the breakdown of time spent in the {\ptop} and {\mtol} kernels for each iteration: although the first iteration is unbalanced, with far field taking about 8 times as much \cpu\ time as near field, later iterations are close to balanced and the total time to solution is optimal. 
We ran extensive tests on the minimum value of $p$ allowed in the relaxed solver, without degrading convergence and accuracy. Table \ref{tab:stokes_min_p} presents some of the data from these tests: for finer surface discretizations, the error degrades when the relaxed value of $p$ is allowed to drop below 3 or 4. To be conservative and avoid any degradation in accuracy, we used a $\pmin=5$ for all further cases with the Stokes equation.

\begin{table}[ht]
\footnotesize
\begin{center}
\begin{tabular}{c|cc|cc|cc}
  & \multicolumn{2}{c|}{2048 panels} & \multicolumn{2}{c|}{8192 panels} & \multicolumn{2}{c}{32768 panels} \\
 $\pmin$ & Error & $it$ & Error & $it$ & Error & $it$ \\ \hline
  & & & & & & \\
 5 & $4.61\times 10^{-2}$ & 24 & $2.48\times 10^{-2}$ & 29 & $1.34\times 10^{-2}$ & 28 \\
 4 & $4.56\times 10^{-2}$ & 25 & $2.50\times 10^{-2}$ & 29 & $1.17\times 10^{-2}$ & 28 \\
 3 & $4.58\times 10^{-2}$ & 23 & $3.57\times 10^{-2}$ & 30 & $1.74\times 10^{-2}$ & 30 \\
 2 & $4.41\times 10^{-2}$ & 25 & $1.45\times 10^{-1}$ & 39 & $2.04\times 10^{-2}$ & 39 \\
 1 & $4.35\times 10^{-2}$ & 28 & $3.42\times 10^{-2}$ & 100\footnotemark[1] & $3.03\times 10^{-2}$ & 100\footnotemark[1] 
\end{tabular}
\end{center}
\caption{Effect of $p_{\text{min}}$ on accuracy and convergence for Stokes flow around a sphere for differing values of $N$. Error is on the total drag force in the $x$-direction, $F_x$.}
\label{tab:stokes_min_p}
\end{table}%

\footnotetext[1]{When $p_{\text{min}}$ is set to 1, solution does not converge within 100 iterations for the two finer meshes.}

\begin{figure}
\begin{center}
	\includegraphics[natwidth=3in,natheight=2in,width=0.5\textwidth]{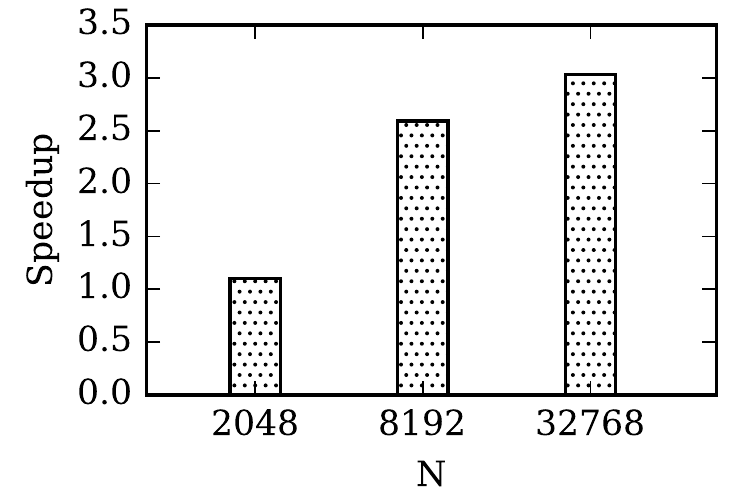}
	\caption{Speed-up for solving first-kind Stokes equation on the surface of a sphere, varying $N$. $10^{-5}$ solver tolerance, $p=16$. Time is measured by averaging the solving time of three identical runs. Plotting script and figure available under CC-BY \cite{WangLaytonBarba2016-figshare3}.}
	\label{fig:stokes_speedup}
\end{center}
\end{figure}

\begin{figure}
	\centering
	\includegraphics[natwidth=3in,natheight=2in,width=0.5\textwidth]{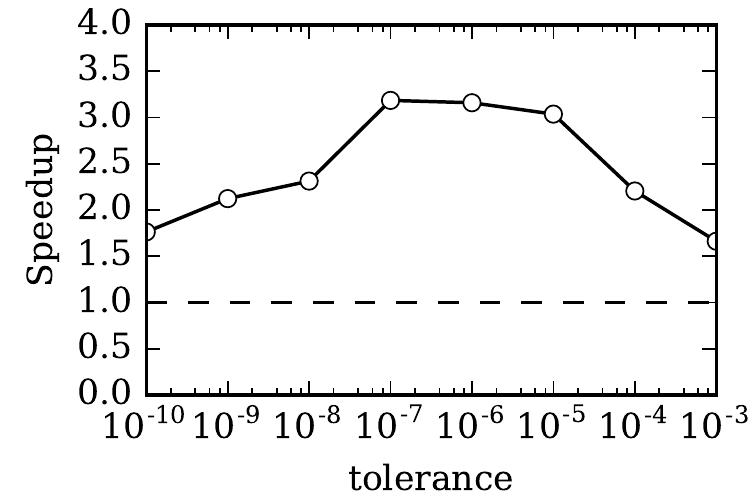}
	\caption{Speed-ups for solving a 1st-kind Stokes problem on a sphere discretized with $8,192$ panels, as the \gmres solver's tolerance increases; $p=16$ for all cases. Time is measured by averaging the solving time of three identical runs. (Multi-threaded evaluator running on 6 \cpu\ cores.) Plotting script and figure available under CC-BY \cite{WangLaytonBarba2016-figshare3}.}
	\label{fig:stokes_tolerance_speedup}
\end{figure}

Figure \ref{fig:stokes_speedup} shows the speed-up resulting from the relaxed \gmres iterations for three increasingly finer surface discretizations. For $N=32,768$, speed-up is more than $3.0\times$. We also observed a speedup between $1.7\times$ and $3.5\times$ under different \gmres solver tolerances for $N=8,192$, as shown in Figure \ref{fig:stokes_tolerance_speedup}.

\subsection{Application to red blood cells in Stokes flow}

A number of medical applications will benefit from greater understanding of the microflows around red blood cells and of the mechanical effects on the cells from this flow. 
The most notable example is perhaps the deadly malaria infection, which makes red blood cells stiffer thus disrupting the flow of blood in capillaries \cite{FedosovETal2011}.
Any design of a biomedical device that processes blood at the micrometer-scale needs to consider the mechanical behavior of blood at the cellular level \cite{Freund2014}. Blood is a dense suspension of mostly red blood cells and smaller concentrations of white blood cells and platelets. The flow regime in small capillaries is at very low Reynolds numbers, and thus completely dominated by viscous effects. 
Red blood cells are very flexible, so any physiologically realistic simulation should take into account their elastic deformations. But here we are only attempting to show the benefit of our relaxation strategy on the Stokes solver, and thus limit our study to the steady Stokes flow around a red-blood-cell geometry. The unsteady problem of coupled Stokes flow and linear elasticity can be approached by repeated solution of boundary-integral problems at every time step, and would equally benefit from the speed-ups seen on a single Stokes solution.

\begin{figure}
\begin{center}
	\subfloat[512 panels]{\includegraphics[natwidth=2.94in,natheight=1.94in,width=0.3\textwidth]{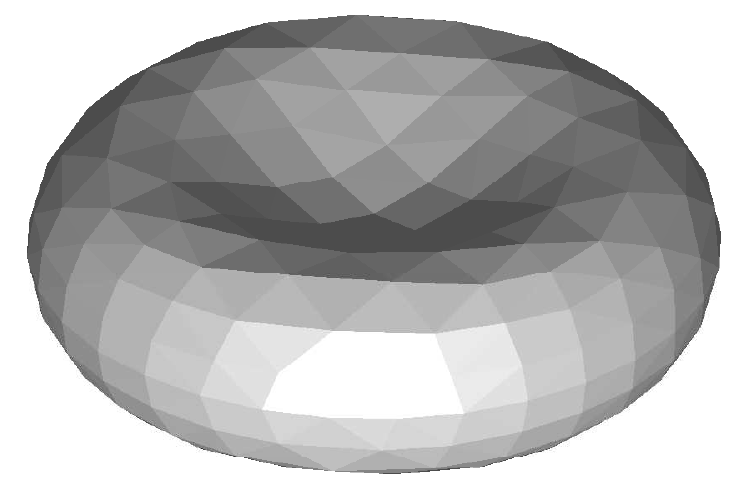}\label{fig:rbc512}}\qquad
	\subfloat[2048 panels]{\includegraphics[natwidth=2.94in,natheight=1.94in,width=0.3\textwidth]{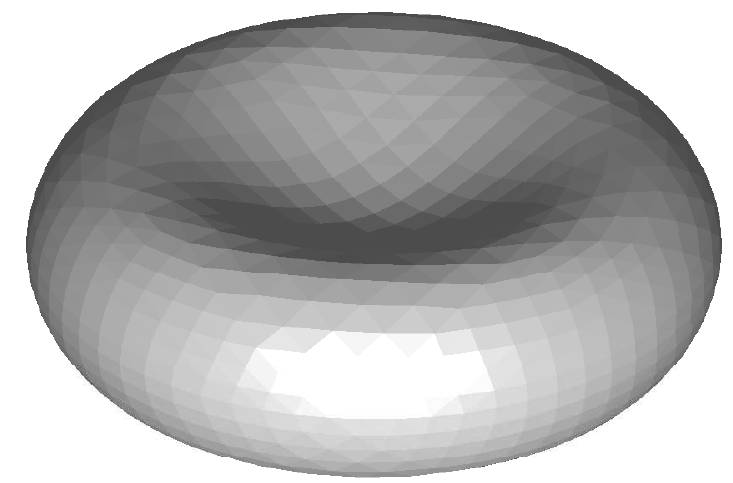}\label{fig:rbc2048}}
	\caption{Surface geometries of red blood cells obtained from transforming sphere triangulations using Equation \eqref{eqn:rbc_parameterization}.}
	\label{fig:glob_rbc}
\end{center}
\end{figure}

To create a surface discretization for a red blood cell, we start with a sphere discretized into triangular panels and transform every vertex $v = v(x,y,z)$, with $x,y,z\in [-1,1]$, into $v' = v'(x',y',z(\rho'))$ using the formula presented in Ref.~\cite{EvansFung1972}:
\begin{equation}
	\label{eqn:rbc_parameterization}
	z(\rho) = \pm \frac{1}{2}\sqrt{1 - \left(\frac{\rho}{r}\right)^{2}}\left ( C_0 + C_2 \left(\frac{\rho}{r}\right)^{2} + C_4\left(\frac{\rho}{r}\right)^{4}\right ),
\end{equation}

\noindent where $x' = x\cdot r,\; y' = y\cdot r,\; \rho = \sqrt{x'^{2}+y'^{2}}$, and the coefficients are: $r=3.91\mu$m,  $C_0= 0.81\mu$m, $C_2= 7.83\mu$m and $C_4=-4.39\mu$m.
Figure \ref{fig:glob_rbc} shows two examples of transformed shapes obtained from sphere triangulations using Equation \eqref{eqn:rbc_parameterization}.

A grid-convergence study using the geometry of a red blood cell requires that we use Richardson extrapolation \cite{roache1998}, since we don't have an analytical solution for this situation. We calculated the drag on a red blood cell in uniform Stokes flow using three surface meshes, consecutively refined by a constant factor $c=4$. 
If the value $f_1$ corresponds to that obtained using the coarsest mesh and $f_2$ and $f_3$ to those using consecutively refined ones, then we can obtain the extrapolated value approximating the exact solution with the following formula:

\begin{equation}
	\bar{f} = \frac{f_1f_3-f_2^{2}}{f_1 -2f_2+f_3}
\end{equation}

\begin{figure}
\begin{center}
	\includegraphics[natwidth=3in,natheight=2in,width=0.5\textwidth]{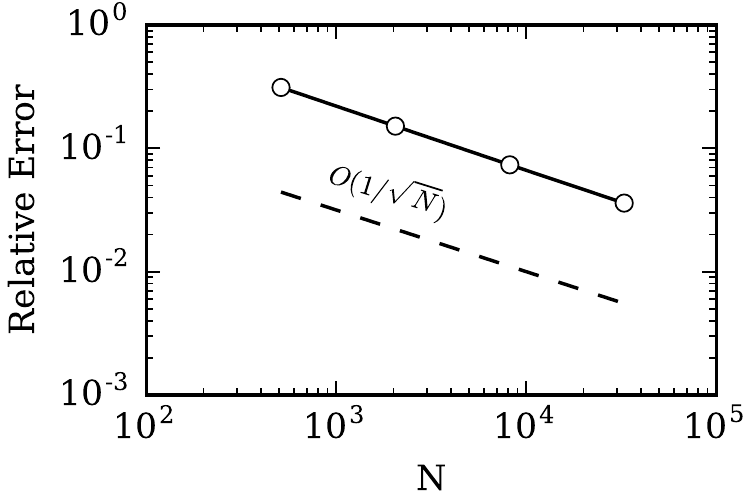}
	\caption{Observed convergence for Stokes flow around red blood cells, with respect to the extrapolated value of the drag coefficient, using Richardson extrapolation \cite{roache1998}. Plotting script and figure available under CC-BY \cite{WangLaytonBarba2016-figshare4}.}
	\label{fig:rbc_extrapolated_convergence}
\end{center}
\end{figure}

Table \ref{tab:rbc_richardson_values} presents the computed values of the drag force obtained with three different meshes, of sizes $N=512$, $2048$ and $8192$. We can also obtain the \emph{observed order of convergence}, $p$, as follows
\begin{equation}
	p = \frac{\ln{\left(\frac{f_2-f_1}{f_3-f_2}\right)}}{\ln{c}},
\end{equation}

\noindent where $c$ is the refinement ratio between two consecutive meshes. With the values in Table \ref{tab:rbc_richardson_values}, the observed order of convergence comes out at $0.52$, matching our expected rate of convergence of $\O{\sqrt{N}}$. 
Figure \ref{fig:rbc_extrapolated_convergence} shows a plot of the error with respect to the extrapolated value, as a function of the mesh size. The plot includes the error obtained with four meshes, with the extrapolated value obtained using the first three coarser meshes.

\begin{table}[h]
\footnotesize
\begin{center}
\begin{tabular}{c|c}
	$N$ & $f_x$ \\
	\hline
	& \\
	$512$ & $-0.057$ \\
	$2048$ & $-0.070$ \\ 
	$8192$ & $-0.077$ \\
\end{tabular}
\end{center}
\caption{Surface mesh sizes and calculated drag force for the convergence study using a red blood cell in uniform Stokes flow.}
\label{tab:rbc_richardson_values}
\end{table}%

\begin{figure}
\begin{center}
	\includegraphics[natwidth=3in,natheight=2in,width=0.5\textwidth]{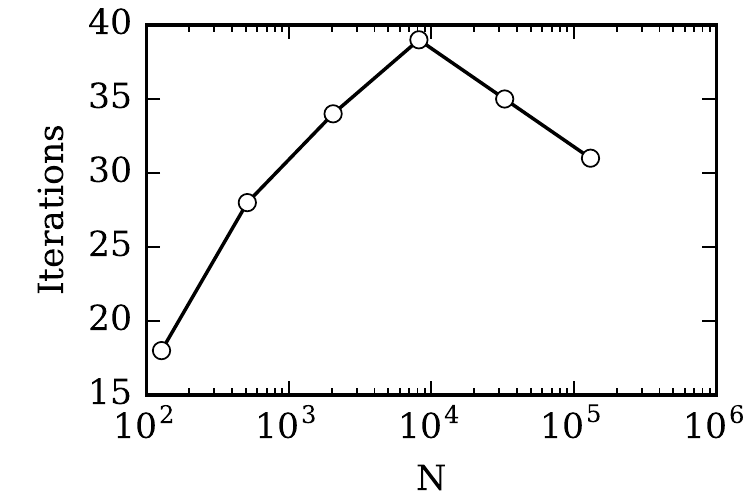}
	\caption{Number of iterations needed for the system to converge for increasingly refined surface meshes on one red blood cell; $p = 16$, target residual $10^{-5}$. Plotting script and figure available under CC-BY \cite{WangLaytonBarba2016-figshare4}.}
	\label{fig:single_cell_iterations}
\end{center}
\end{figure}

The calculations with increasingly finer surface meshes take more time to complete not only because the number of unknowns is larger, but also because they may require a greater number of iterations to converge to a desired residual.
Figure \ref{fig:single_cell_iterations} shows that the number of iterations needed as the surface mesh varies from size $N=128$ to $8,192$ increases from 18 to 39. For further refined meshes, the number of iterations remains about the same. This is an indication that the surface mesh is sufficiently refined with $8,192$ panels for one red blood cell.

\begin{figure}[ht]
\begin{center}
\includegraphics[natwidth=6in,natheight=2.88in,width=0.75\textwidth]{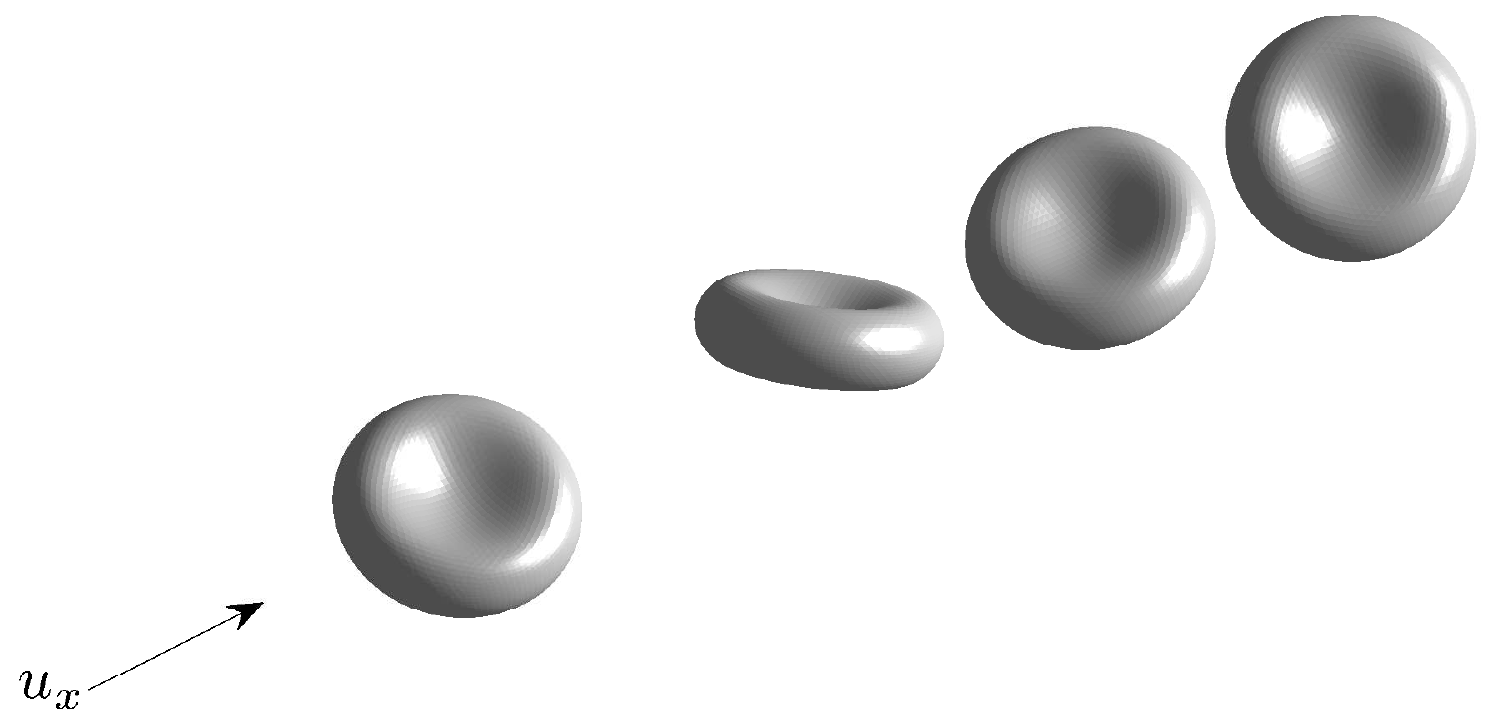}
	\caption{Surfaces for four red blood cells in a uniform Stokes flow.}
	\label{fig:multiple_cells}
\end{center}
\end{figure}

For the next tests, we used several red blood cells in a sparse spatial arrangement. Realistic blood flows have densely packed red blood cells, but the purpose of this test is to simply demonstrate the boundary element solver with a larger problem. We set up a collection of red blood cells by making copies of a discretized cell, then randomly rotating each one, and shifting it spatially in each  coordinate direction by a positive random amount that ensures they do not overlap. The resulting arrangement may look like that shown in Figure \ref{fig:multiple_cells}.

\begin{figure}
\begin{center}
	\includegraphics[natwidth=7in,natheight=3.5in,width=0.85\textwidth]{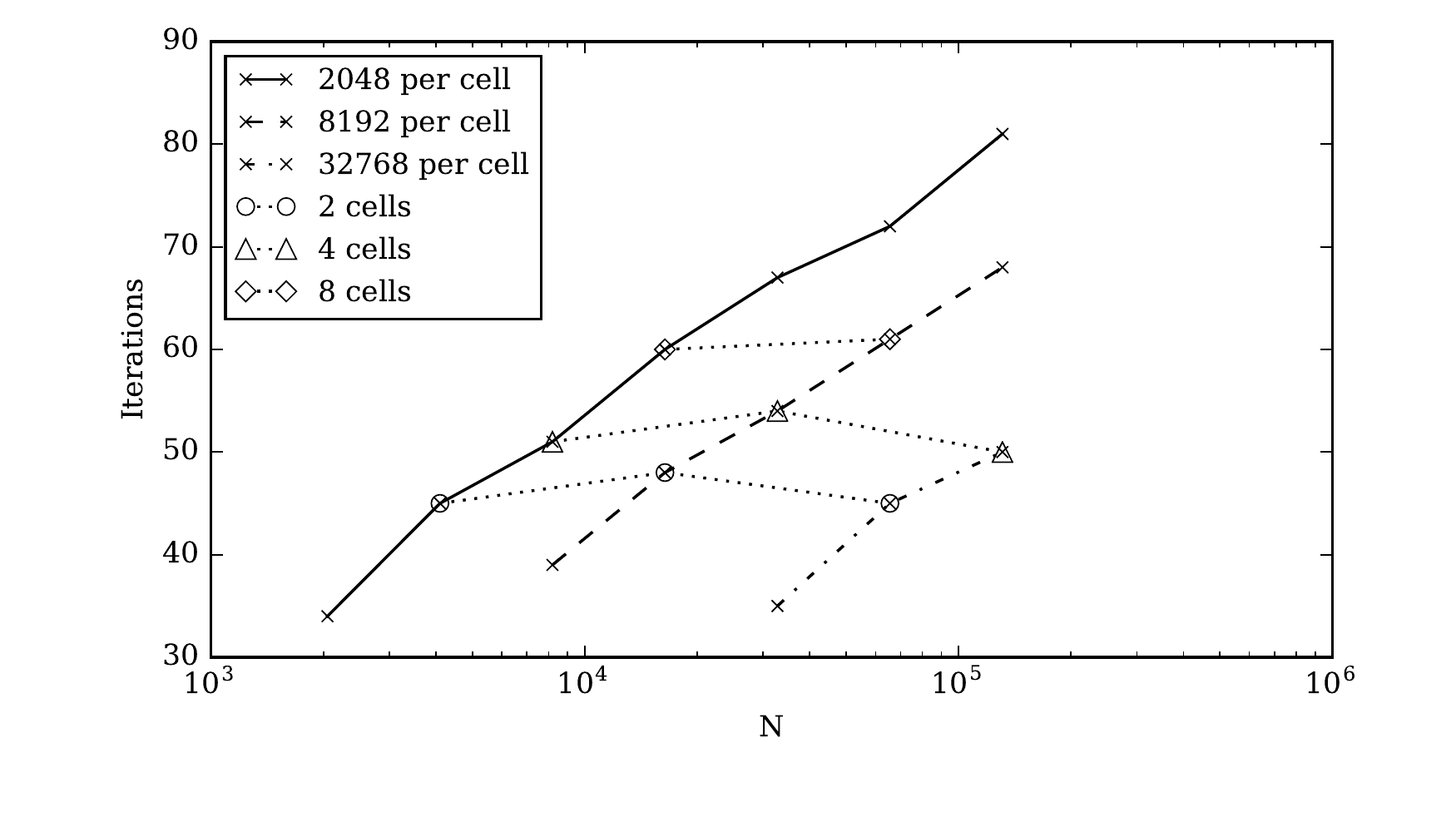}
	\caption{Number of iterations needed to converge to a desired residual of $10^{-5}$ for systems with multiple red blood cells, discretized with different mesh sizes ($p = 16$). Plotting script and figure available under CC-BY \cite{WangLaytonBarba2016-figshare4}.}
	\label{fig:multiple_cell_iterations}
\end{center}
\end{figure}

Using 2, 4 and 8 red blood cells, we looked at the number of iterations needed to converge to a solver tolerance of $10^{-5}$ when using three different surface mesh sizes on each cell: $N=2048$, $8192$ and $32,768$. In all cases $p = 16,\;p_{\text{min}} = 5$. Figure \ref{fig:multiple_cell_iterations} shows that the number of iterations needed to converge increases sharply with the number of red blood cells in the system, while the number of panels per cell has a smaller effect in this range of mesh sizes. In all cases, the number of iterations is between 45 and 65, and thus we expect to see good speed-ups using the relaxation strategy.

\begin{figure}[ht]
\begin{center}
	\includegraphics[natwidth=4in,natheight=3in,width=0.55\textwidth]{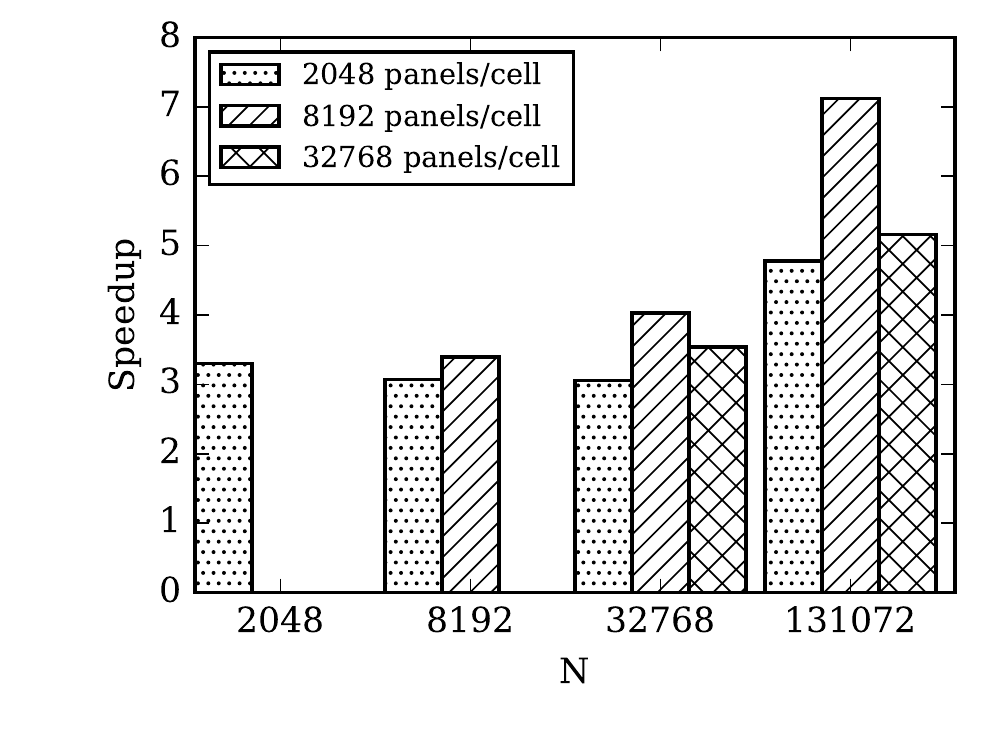}
	\caption{Speed-ups of the relaxation strategy with several red blood cells in uniform Stokes flow, using different mesh sizes on each cell. The abscissa value corresponds to the total number of panels (all cells). Test parameters shown in Table \ref{tab:cells_relaxation_settings}. Time is measured by averaging the solving time of three identical runs.} 
	\label{fig:multiple_cell_speedup}
\end{center}
\end{figure}

We completed several tests of the relaxation strategy, using between 1 and 64 red blood cells, and various mesh sizes on each cell. The common parameter settings for these tests are listed in Table \ref{tab:cells_relaxation_settings} and, in each case, we chose the value of $\ncrit$ (establishing the balance between near and far field) to obtain the smallest time to solution for that run. The detailed results are listed in Tables \ref{tab:single_cell_relaxation_results}, \ref{tab:multiple_cell_relaxation_results_2048}, \ref{tab:multiple_cell_relaxation_results_8192} and \ref{tab:multiple_cell_relaxation_results_32768}, and Figure \ref{fig:multiple_cell_speedup} shows a summary of the observed speed-ups, mostly hovering close to $4\times$. The largest problems, with a total of $131,072$ panels (all cells combined), are atypical because in these cases we are unable to use an efficient sparse-matrix representation of the near field, due to the large memory requirement. But this can also be seen as an advantage of the relaxation strategy, which leads to using smaller near fields and thus extends the range of problem sizes where we can use the efficient sparse-matrix representation. Indeed, if one needed to solve a problem of size $N=131,072$ (on one \cpu\ using six threads, like we do here), then the potential for a $7\times$ speed-up is real.

\begin{table}[h]
\footnotesize
\begin{center}
\begin{tabular}{c|c}
 Variable & Setting \\ 
\hline
 & \\
 $p_{\text{initial}}$ & $16$ \\
 $p_{\text{min}}$ &  $5$ \\
 solver tolerance & $10^{-5}$ \\
 Near-field & Sparse matrix \\
 Threads & $6$ \\
 Solver & {\gmres} \\ 
 Preconditioner & None
\end{tabular}
\end{center}
\caption{Parameters for the tests of the relaxation strategy with red blood cells in uniform Stokes flow.}
\label{tab:cells_relaxation_settings}
\end{table}%

\begin{table}[htp]
\footnotesize
\begin{center}
\begin{tabular}{c|c|c|c|c|c|c}
 & & \multicolumn{2}{c|}{Non-relaxed} & \multicolumn{2}{c|}{Relaxed} \\
 & & \multicolumn{2}{c|}{} & \multicolumn{2}{c|}{} \\
 N & \# unknowns & $\ncrit$ & $t_{\text{solve}}$ & $\ncrit$ & $t_{\text{solve}}$ & Speedup \\\hline
 & & & & & \\
 2048 & 6144 & 400 & 48.18 & 150 & 14.61 & 3.30 \\
 8192 & 24576 & 400 & 234.44 & 100 & 69.03 & 3.40 \\
 32768 & 98304 & 400 & 978.80 & 100 & 276.52 & 3.54 \\
 131072 & 393216 & 200 & 5122.43\footnotemark[1] & 100 & 1044.77 & 4.90\footnotemark[1] \\	
\end{tabular}
\end{center}
\caption{Timings and speed-up of the relaxation strategy on a single red blood cell in uniform Stokes flow, with the test parameters shown in Table \ref{tab:cells_relaxation_settings}.}
\label{tab:single_cell_relaxation_results}
\end{table}%

\footnotetext[1]{Due to memory restrictions, a sparse-matrix representation of the near-field could not be used, resulting in a much slower {\ptop} evaluation.}

\begin{table}[htp]
\footnotesize
\begin{center}
\begin{tabular}{c|c|c|cc|cc|c}
\multicolumn{8}{c}{2048 panels / cell} \\
& & & \multicolumn{2}{c}{Non-relaxed} & \multicolumn{2}{c}{Relaxed}\\
N & \# unknowns & $N_c$ & $\ncrit$ & $\tsolve$ & $\ncrit$ & $\tsolve$ & Speedup \\ \hline
& & & & & & &  \\
2048 & 6144 & 1 & 400 & 48.18 & 150 & 14.61 & 3.30 \\ 
8192 & 24576 & 4 & 400 & 267.64 & 100 & 87.20 & 3.07 \\ 
32768 & 98304 & 16 & 400 & 1358.07 & 100 & 444.55 & 3.05 \\
131072 & 393216 & 64 & 200 & 7077.77\footnotemark[1] & 100 & 1481.47 & 4.78\footnotemark[1] \\
\end{tabular}
\end{center}
\caption{Timings and speed-up of the relaxation strategy with several red blood cells in uniform Stokes flow, each cell discretized with 2048 panels and test parameters shown in Table \ref{tab:cells_relaxation_settings}.}
\label{tab:multiple_cell_relaxation_results_2048}
\end{table}

\begin{table}[htp]
\footnotesize
\begin{center}
\begin{tabular}{c|c|c|cc|cc|c}
\multicolumn{8}{c}{8192 panels / cell} \\
& & &  \multicolumn{2}{c}{Non-relaxed} & \multicolumn{2}{c}{Relaxed}\\
N & \# unknowns & $N_c$ & $\ncrit$ & $\tsolve$ & $\ncrit$ & $\tsolve$ & speedup \\ \hline
& & & & & & &  \\
8192 & 24576 & 1 & 400 & 234.44 & 100 & 69.03 & 3.40 \\ 
32768 & 98304 & 4 & 400 & 1494.37 & 100 & 371.06 & 4.03 \\
131072 & 393216 & 16 & 150 & 12817.01\footnotemark[1] & 100 & 1800.53 & 7.12\footnotemark[1] \\
\end{tabular}
\end{center}
\caption{Timings and speed-up of the relaxation strategy with several red blood cells in uniform Stokes flow, each cell discretized with $8,192$ panels and test parameters shown in Table \ref{tab:cells_relaxation_settings}.}
\label{tab:multiple_cell_relaxation_results_8192}
\end{table}

\begin{table}[htp]
\footnotesize
\begin{center}
\begin{tabular}{c|c|c|cc|cc|c}
\multicolumn{8}{c}{32768 panels / cell} \\
& & & \multicolumn{2}{c}{Non-relaxed} & \multicolumn{2}{c}{Relaxed}\\
N & \# unknowns & $N_c$ & $\ncrit$ & $\tsolve$ & $\ncrit$ & $\tsolve$ & speedup \\ \hline
& & & & & & &  \\
32768 & 98304 & 1 & 400 & 978.80 & 100 & 276.52 & 3.54 \\
131072 & 393216 & 4 & 200 & 7209.57\footnotemark[1] & 100 & 1397.27 & 5.16\footnotemark[1] \\	
\end{tabular}
\end{center}
\caption{Timings and speed-up of the relaxation strategy with several red blood cells in uniform Stokes flow, each cell discretized with $32,768$ panels and test parameters shown in Table \ref{tab:cells_relaxation_settings}.}
\label{tab:multiple_cell_relaxation_results_32768}
\end{table}

\subsection{Reproducibility}

Our inexact \gmres code is open-source and available on Github. Knowing the increasing importance of computational reproducibility, we provided the ``\textit{reproducibility package}" --- containing the running and post-processing scripts --- to generate the figures in this result section. To save readers from potential headaches of dependency mismatches, we also prepared a Dockerfile to setup up the software environment (a Docker container), under which the ``\textit{reproducibility package}" can be used to replicate our results.

\section{Conclusion} 

We have shown the first successful application of a relaxation strategy for fast-multipole-accelerated boundary element methods, based on the theory of inexact \gmres. Testing the relaxation strategy on Laplace problems, we confirmed that it converges to the right solution, it provides moderate speed-ups over using a fixed $p$, and it leads to initially bloated far-fields to obtain the minimum time to solution.
Exploring the performance advantage of relaxing the value of $p$ as \gmres iterations advance, we concluded that problems requiring high accuracy and/or resulting in more ill-conditioned linear systems will experience the best speed-ups, which for Laplace problems were between $2.1\times$ and $3.3\times$ in our tests on a sphere with constant potential.

In the case of the Stokes equation, the speed-ups that can be obtained using a relaxation strategy are larger, due to the fact that Stokes problems require both more iterations to converge (and the relaxed solver spends more time at low $p$) and more work per iteration (equivalent to four Laplace evaluations). 
We found that it's important for Stokes problems to also enforce a minimum value of $p$ to avoid accuracy or convergence degradation.
Relaxed \gmres iterations in this case reduced the time to solution by up to $3.0\times$ in tests of Stokes flow around a sphere. 
We completed various tests for Stokes flow around red blood cells, with up to 64 cells. We studied numerical convergence in this situation using Richardson extrapolation and obtained an observed order of convergence of $0.52$, close to the expected value of $1/2$. The speed-ups resulting from the relaxation strategy in these tests were in most cases close to $4\times$.

Relaxing the truncation order $p$ in the multipole expansions as Krylov iterations progress is one of those seemingly simple ideas that strike one as obvious a posteriori. Yet, as far as we know, it has not been tried before, nor has it been implemented in a \bem. Given this method's wide popularity in computational engineering, we look forward to many applications benefitting from healthy speed-ups from applying relaxed-$p$ \fmm. We  showed that Stokes problems, in particular, can expect $4\times$ speed-ups in large problems still fitting on one workstation. Linear elasticity problems should experience similar speed-ups (although we didn't try them). This is pure algorithmic speed-up that should multiply with any hardware speed-ups obtained, for example, by moving computational kernels to \gpu s.

\appendix
\section{Algorithm listings}\label{sec:algorithms}
 
 \begin{algorithm}
 \footnotesize
	\caption{Matrix-vector multiplication.}
	\label{alg:matvec}
	\begin{algorithmic}
		\State Initialize $\mathbf{w}$
		\For{Collocation points $i=1\cdots N_p$}
			\State $w_i \gets 0$
			\For{Integration panels $j=1 \cdots N_p$}
				\For{Gauss quadrature points $k= 1 \cdots N_k$}
				\State $w_i \gets w_i + v_j \cdot q_k \cdot S_j \, \frac{1}{|\mathbf{x_i}-\mathbf{x_k}|}$
				\EndFor
			\EndFor
		\EndFor 
	\end{algorithmic}
\end{algorithm}

\clearpage
\section*{Acknowledgements}
 This work was supported by the National Science Foundation via NSF CAREER award OCI-1149784. LAB also acknowledges support from NVIDIA Corp.\ via the CUDA Fellows Program. 
Dr. Cris Cecka (previously at Harvard University, currently at Nvidia Corp.) contributed to the development of the code, particularly writing the octree and base evaluator. He later continued developing his fast-multipole framework, which evolved into his FMMTL project; see \href{https://github.com/ccecka/fmmtl}{https://github.com/ccecka/fmmtl}.
The authors also wish to acknowledge valuable interactions with Dr. Christopher Cooper (previously at Boston University, currently at Universidad T{\'e}cnica Federico Santa Mar{\'i}a) that helped with the implementation of the boundary element method.

\section*{References}
\bibliographystyle{elsarticle-num}

\end{document}